\documentclass[dvips]{imsart}
\makeindex
\usepackage{amsthm,amsmath}

\startlocaldefs
\usepackage{latexsym,amsmath,amsfonts}
\usepackage{mathrsfs}
\newcommand{\nref}[1]{(\ref{#1})}
\def\Cal{\cal}
\renewcommand{\textbf}[1]{\begingroup\bfseries\mathversion{bold}#1\endgroup}

\def\text#1{\hbox{#1}}

\newcommand{\R}{\mathbb R}

\newcommand{\E}{\mathbb E}
\def\1{\mbox{1\hspace{-.20em}I}}
\newtheorem{theorem}{Theorem}[section]
\newtheorem{proposition}{Proposition}[section]
\newtheorem{lemma}{Lemma}[section]

\newtheorem{remark}{Remark}[section]

\newcommand{\CN}{{\Cal{N}}}

\def\a{\alpha}

\def\e{\epsilon}
\def\r_\e{r_\epsilon}

\def\g{\gamma}

\def\l{\left}
\def\r{\right}

\def\Var{{\rm Var}}

\numberwithin{equation}{section}
\def\AArm{\fam0 \mathrm}%
\def\AAk#1#2{\setbox\AAbo=\hbox{#2}\AAdi=\wd\AAbo\kern#1\AAdi{}}%
\def\AAr#1#2#3{\setbox\AAbo=\hbox{#2}\AAdi=\ht\AAbo\raise#1\AAdi\hbox{#3}}%
\def\BBn{{\AArm I\!N}}%
\def\BBr{{\AArm I\!\!R}}%
\def\BBz{{\AArm Z\!\!\!Z}}%
\def \E{\hbox{\it I\hskip -2pt E}}
\def \P{\hbox{\it I\hskip -2pt P}}
\def \V{\mbox{\rm Var}}
\def \I{\hbox{\rm 1\hskip -3pt I}}

\endlocaldefs
\begin{document}

\begin{frontmatter}

\title{Detection of sparse  additive functions}
\runtitle{Detection of sparse functional signals}

\begin{aug}
\author{\fnms{Ghislaine} \snm{Gayraud}\thanksref{t1} \corref{Gayraud}\ead[label=e1]{ghislaine.gayraud@utc.fr}}

\address{Universit\'e de Technologie de Compi\`egne $\&$ CREST, \\  BP 20529, 60205 Compi\`egne, France\\ email: ghislaine.gayraud@utc.fr}

\and
\author{\fnms{Yuri} \snm{Ingster}\thanksref{t2}\ead[label=e2]{yurii_ingster@mail.ru}}
\address{St. Petersburg State Electrotechnical University, 5, \\ Prof. Popov str., 197376 St.Petersburg, Russia \\ email:
yurii$\_$ingster@mail.ru }

\thankstext{t1}{Ghislaine Gayraud's research was partially supported by the
ANR-blanc 'SP Bayes'}
\thankstext{t2}{Yuri Ingster's research was partially supported  by RFBR Grant
11-01-00577 and by Grant NSh--4472.2010.1}

\runauthor{G. Gayraud  et al.}

\affiliation{Universit\'e de Technologie de Compi\`egne $\&$ CREST \thanksref{t1}\\ Petersburg State Electrotechnical University \thanksref{t2}}

\end{aug}





\begin{abstract}
 We study  the problem of detection of   high-dimensional
signal functions in the Gaussian white  noise model. We assume  that,
in addition to
a smoothness assumption, the signal function has an additive sparse structure. The detection problem is expressed in terms of
a nonparametric hypothesis testing problem and is solved
using asymptotically minimax approach.  We provide  minimax
test procedures that are adaptive in the sparsity parameter in the
high sparsity case. We extend some known
results related to the detection of sparse high-dimensional vectors to the functional case. In
particular, our derivation of asymptotic detection rates is based on
same detection boundaries as in the vector case.
\end{abstract}




\tableofcontents

\end{frontmatter}

\section{Introduction}\label{Intro}
Over the past years, boosted by applications and computer performance,  problems in high-dimensions have been explored in a number of statistical studies. If no additional structure is assumed, high-dimensional data processing  suffers from some intrinsic difficulties such as the
curse of dimensionality that results in  a  loss in the efficiency of statistical procedures, and inconsistency of classical statistical procedures -- even in the
linear regression model -- unless the dimension of variables is less than the sample size.

In order to overcome the curse of dimensionality in a nonparametric framework, where typical functional classes are  Sobolev, Holder,  or Besov balls,
some additional conditions, including additivity or tensor product
structure, are
assumed, see, for instance, \cite{{St.85}, {IH.97}, {L.00},
{IS.05}, {IS.07a}, {IS.07b}} and references therein.
Even if  one of these conditions is assumed, yet it is required that the sample size is to be larger than the data dimension.
One way to free oneself from the latter condition is to impose an additional sparsity constraint. \\

In this paper we focus on the problem of  detection of high-dimensional signal functions in the Gaussian
white  noise model.  To avoid   difficulties stemming from  high-dimensional  settings,
we suppose that a signal function satisfies an additional structural condition. Specifically, it is assumed to
be sparse additive.
This means that a high-dimensional
function of interest is a  sum of few univariate functions.
Formally, we consider an $d$-dimensional ($d \in \BBn$ and $d >0$) Gaussian white  noise model
\begin{eqnarray}\label{M1}
dX(t)=f(t)dt+\epsilon dW(t),\ t\in [0,1]^d,
\end{eqnarray}
where $W(t)$ is the Wiener process, $\epsilon>0$ is the noise level,
and $f$, the quantity of interest,  is the signal function. The additive sparse structure
 means that $f$  is the sum of $d$ univariate functions $f_j$:
\begin{eqnarray}
f(t)=\sum_{j=1}^d \xi_jf_j(t_j),\quad t_j \in [0,1], \label{add-struc}
\end{eqnarray}
where the $\xi_j$'s are unknown but deterministic taking their values in $\{0,1\}$ : ``0''  means that   the $j$th
component $f_j$ is non active  whereas  ``1'' means that $f_j$ is active.
Denote by $K$ the positive number of active components, that is, $K=\displaystyle{\sum_{j=1}^d} \xi_j$,  and assume that $K=d^{1-b}$,
where $b\in (0,1)$ is the {\it sparsity index}. If $d^{1-b}$ is not an integer then take $K$ as its integer part.
Denote by ${\cal F}_{d,b}$ the functional class of additive sparse signals $f$ of the
form  (\ref{add-struc}) with $K=d^{1-b}$ active components and $d^b$ non-active components.
Model (\ref{M1}) with the sparse additive structure (\ref{add-struc}) is a natural generalization of the sparse linear model: the nonparametric
nature of the problem suggests to consider more flexible models.

 There is a huge statistical literature on estimation in  sparse  models, see, for
instance, \cite{{BRT.09},{KJL.07},{D.06}} and references therein. In  particular, there are  many works related to the well-known Lasso introduced  by Tibshirani  \cite{T.96} in 1996.
There are also a number of papers that deal with nonparametric estimation in sparse additive models.
For a complete review of these topics, we refer to  \cite{RWY.11}, where minimax estimation rates in sparse additive models are obtained, to \cite{HHF.10},
where the Lasso-type estimate in sparse additive models is studied, and to \cite{St.85}, where various structural assumptions on
models in high dimensions are discussed. \\

Back to our study, the detection problem at hand can be expressed in terms of a nonparametric hypothesis testing problem with the null hypothesis
stating that ``the signal is a constant'', and ``there is no signal'' being a particular case of the null hypothesis.
In
order to specify an alternative hypothesis, recall that, within the minimax framework, it is impossible
to detect  signal functions that are ``too close'' to the null one, as well as to test the null and alternative hypotheses for too large  alternative classes.
Therefore, we are interested in the following nonparametric hypothesis testing problem:
\begin{eqnarray}H_0\; :\; f= const_0 \;\;\;\;\; \mbox{{\rm
versus }} \;\;\;\;\;  H_{1} \;: \; f=const_1+f^{1}, \; f^{1} \in {\cal
F}_d(\tau, r_\epsilon, b),\label{PB-Test}\end{eqnarray}
where
$$\left\{
\begin{array}{l}
const_0, \; const_1 \mbox{ {\rm are some constants,} } \\
 {\cal
F}_d(\tau, r_\epsilon, b)= \!\! \{f \in {\cal F}_{d,b} :\forall j,  f_j \in \tilde S_\tau   \mbox{ {\rm and } } \|f_j\|_2
\geq r_\epsilon \}, \; \tau >0, \; r_\epsilon >0, \\
\tilde S_\tau= \{ f \in L_2([0,1]): \; \int_0^1 f(t)dt=0, \;
\|f\|_2^{(\tau)}\le 1\}. 
\end{array}\right.
$$
The $L_2$-norm $\|\cdot\|_2$ is used to separate the nonparametric alternative from the null hypothesis.
The functional class  $\tilde S_\tau$
 is the Sobolev ball, expressed via the  Sobolev semi-norm $\|\cdot\|_2^{(\tau)}$,  that contains  $\tau$-smooth functions, which are
assumed 1-periodic and
orthogonal to a constant. Due to the periodic constraints, it is possible to express $\|\cdot\|_2^{(\tau)}$   in terms of Fourier coefficients; this  will be done in Section \ref{sec:Bgk-Seq-Sp}.
 The quantity  $\tau$   is the smoothness parameter.
Both the smoothness condition and the  separation condition between $H_0$ and $H_1$ are expressed in terms of  the components $f_j$ that are linked to the whole signal $f$ via
(\ref{add-struc}): each active component $f_j$ is smooth and is separated from the null hypothesis in the $L_2$-norm by a positive value
$r_\epsilon$.  

In Section  \ref{sec:Extended}, we  generalize the  hypothesis testing problem (\ref{PB-Test})  by considering a more  general class of alternatives that
consists of signals $f$ equal,  up to a constant, to  a function $f^{1} \in {\cal F}_{d,b}$, which is separated from the null hypothesis
 in the $L_2 ([0,1]^K)$-norm, and whose smoothness is expressed in terms of the whole function $f$.

%
%
For these two hypothesis testing problems, the main questions are:
{\it what are the separation rates in the problem, i.e., what are
the asymptotics for the minimal $r_\epsilon$ such that one can
distinguish between $H_0$ and $H_1$?  And, also, what are the optimal test
procedures that provide distinguishability}?

To answer these questions, we use asymptotically minimax approach that provides detection boundaries or distinguishability conditions, i.e.,  necessary and sufficient conditions for the possibility of successful detection; these detection boundaries
yield
 asymptotics for the minimal  $r_\e$ separating
the areas of  distinguishability  and non-distinguishability (between $H_0$ and $H_1$).
The asymptotics for the minimal  values of $r_\e$  are called  either the (minimax) separation rates or the minimax rates of testing; in the present paper,  the separation rates are denoted  by $r_\e^\star$.

In connection with the current study,  a number of  works   on detection and classification boundaries
in Gaussian sequence models could be mentioned,  see, for example,  \cite{{I.97},
{I.01}, {I.02}, {IS.02a}, {IS.02b}, {DJ.04}, {IS.07a},{IS.07b}, {IPT.09}}.
 Also, in \cite{ITV.10}, rather than considering a Gaussian sequence model,  the authors generalize the problem of finding a detection boundary in the linear regression model.
Another paper \cite{IL.03} deals with the signal detection problem in a multichannel model in the functional framework. At the end of the next paragraph,
we explain what are the  differences between the results in \cite{IL.03} and our study. \\

The main contribution  of this paper
consists of extending the results on detection boundaries obtained for $d$-dimensional sparse Gaussian vectors, see, for instance, \cite{IS.02b}, to the functional case.
In particular, we obtain the same detection  boundaries  as in the vectorial case.
 However, in the case of high sparsity when $b>1/2$,
 an additional assumption on the growth of  $d$ as a function of $\epsilon$ is required.
Distinguishability is possible when the sum of the  type I error probability and the maximum  over  alternatives of the type II error probability
vanishes asymptotically, and
distinguishability  is not possible when this sum tends  to one.
%
Boundary conditions depend on the quantity $a(r_\epsilon)=a(r_\epsilon,d,\tau)$,
which is a solution
of a certain  extremal problem stated in Section \ref{sec:Extreme}. In the vectorial case, the quantity $a(r_\epsilon)$ corresponds
to the energy of a signal (see \cite{IS.02b} and \cite{IL.03}). In the functional case, this quantity characterizes the distinguishability in a one-variable
hypotheses testing problem.  The minimax separation rates obtained in this paper depend on the value of $b$: for large $b$ they are worse than
for small $b$. Such a behaviour is expected because, with large  $b$, only few components are active, and hence
the problem of distinguishing between the alternative and null hypothesis becomes more difficult.

For the most difficult case of $b\in(1/2,1)$, not only separation
rates, but also sharp separation rates, that include both rates and constants, are obtained. 
We also provide  optimal test procedures for which
 minimax rates of testing are achieved asymptotically.
Depending on the value of $b$, we propose two types of test procedures: one is of a $\chi^2$ type,
the other one is related to  a  Higher-Criticism statistic
introduced in \cite{DJ.04} and based on the Tukey's ideas.
In the case of $b\in(1/2,1),$ our test procedure  is adaptive in the sparsity index  $b$, see Remark \ref{adapt}.

In the paper \cite{IL.03}, which is focused on a similar  problem of multichannel signal
detection, the optimal rates are obtained. In our study,
we obtain sharp separation rates for $b\in (1/2,1)$.
The main difference between the study of \cite{IL.03} and our work is in the
quantity $a(r_\epsilon)$ that characterizes the distinguishability: in our work, it is just a solution of a certain extremal problem,
whereas in \cite{IL.03}, it is obtained directly from the use of
the respective test procedures. 
\\

%
%


The rest of the paper is organized as follows.
Section \ref{sec:Bgk-Seq-Sp} is concerned with the problem of finding detection boundary in a sparse Gaussian $d$-vectors model.
In Section \ref{sec:Trans}, we give  a new  formulation of the problem (\ref{PB-Test}) in terms of sequence spaces.
Section \ref{sec:Extreme} is devoted to
the description of the extremal problem that gives the distinguishability characteristics. 
The main results are stated
in Section \ref{sec:Main}. In
Section \ref{sec:Extended}, we generalize  the hypothesis testing problem (\ref{PB-Test})  by  considering
more general alternatives. The proofs are given in Section \ref{sec:Proofs}.

\section{Detection boundaries in a vectorial Gaussian model} \label{sec:Bgk-Seq-Sp}


Hypothesis testing problems for $d$-dimensional vectors,
under the sparse conditions similar to the ones we use, were studied in
\cite{{I.97},{IS.02b},{DJ.04}}. Namely, let $X=(X_1,\ldots,X_d)$ be a
random vector of the form $X_j=v_j+\eta_j,$ where
$\eta_j \stackrel{\mbox{{\rm i.i.d.}}}{\sim} \CN(0,1)$, $j=1,\ldots,d$, and
\begin{equation}\label{v1}
v_j=\xi_ja,\quad a>0,\quad \xi_j \in\{0,1\},\quad
K=\sum_{j=1}^d\xi_j=d^{1-b},\quad b\in (0,1).\qquad
\end{equation}
Let $V_d(a,b)\subset \R^d$ be the set of all vectors
$v=(v_1,\ldots,v_d)$ of the form \nref{v1}. Then, the testing problem is stated as follows: it is required to test $H_0: v=0$ against the alternative $H_1: v\in V_d(a,b)$.
Here the questions of interest are: what are the asymptotics for
$a=a_d$ as $d\to  + \infty$
for which the hypotheses $H_0$ and $H_1$ separate asymptotically?
Also,  what are the optimal test procedures that provide the  distinguishability (or separation) of $H_0$ and $H_1$?

The answer to each question depends essentially on the sparsity index $b\in (0,1)$,
see \cite{{I.97},{IS.02b},{DJ.04}}.
The detection boundaries are expressed in terms of $a$, $d$ and $b$:  if $b \leq 1/2$ (moderate sparsity),
then the distinguishability is impossible when $ad^{1/2-b}=o(1)$, and it
is possible when $ad^{1/2-b}\to + \infty$.  This is achieved by the test procedure
based on a simple linear statistic $t=d^{-1/2}\sum_{i=1}^dX_i$. If
$b>1/2$ (high sparsity), then the distinguishability conditions look as follows:
the distinguishability is impossible when $\limsup a/T_d<\varphi(b)$,
and it is possible when $\liminf a/T_d>\varphi(b)$, where $T_d=\sqrt{\log(d)}$
and the function $\varphi(b),\ b\in (1/2,1)$ is defined by
\begin{equation}\label{phi}
 \varphi(b)=\begin{cases}\varphi_1(b)=\sqrt{2b-1},& 1/2<b\le 3/4,\\
     \varphi_2(b)=\sqrt{2}(1-\sqrt{1-b}),& 3/4<b<1.
     \end{cases}
\end{equation}
Observe that the function $\varphi$ is positive, continuous, and increasing in
$b\in (0,1]$.

The test procedure that provides distinguishability in the high-sparsity case is based on the
Higher-Criticism statistics introduced in \cite{DJ.04}. It is
defined  as $ L_d=\displaystyle{\max_{s>s_0}} L_d(s) $, for
any $s_0>0$, with
\begin{eqnarray}
L_d(s)=\frac{1}{\sqrt{d\,\Phi(s)\Phi(-s)}}\sum_{i=1}^d(\I_{(X_i>s)}-\Phi(-s)), \label{HC-DJ04}
\end{eqnarray}
where, here and later, $\Phi$ stands for the standard Gaussian cumulative distribution function. Note that it suffices to
take the maximum of $L_d$ over a discrete grid of the form
$s_l=u_lT_d,\ u_l=\delta_d l,\ l=1,\ldots,L$, such that $u_L\leq
\sqrt{2}$ and $\delta_d=o(1)$ is small enough.

\section{Transformation of the statistical testing problem} \label{sec:Trans}
Consider the tensor structure of the space $L_2([0,1]^d)$, i.e., $L_2([0,1]^d)=L_2([0,1]) \otimes \ldots \otimes L_2([0,1])$.  Then,   the corresponding  orthonormal basis
 $({\tilde \phi}_{l}^d)_{l \in \BBz^d}$ of
$ L_2 ([0,1]^d)$ has the form
$${\tilde \phi}_{l}^d (t) =\prod_{j=1}^d \phi_{l_j}^1(t_j), \; t=(t_1,\ldots,t_d) \in [0,1]^d, l=(l_1,\ldots,l_d) \in \BBz^d, $$
 where
 $(\phi^1_k)_{k \in { \BBz}}$ is an orthonormal basis of $L_2([0,1])$.
It is assumed that $\phi^1_0=1$.
For any $(j,k) \in \{1,\ldots,d\} \times \BBz$, let us define ${\bar \phi}_{j,k}^d$ as
$$ {\bar \phi}_{j,k}^d(t)={\tilde \phi}_l^d(t)=\phi_k^1 (t_j), \;l=(0,\ldots,k,0,\ldots,0),$$ where $k$ is the $j$-th component of $l$. Observe that
$ {\bar \phi}_{j,0}^d=1$.
Using the orthonormal system $({\bar \phi}_{j,k}^d)_{(j,k)\in \{1,\ldots,d\} \times \BBz}$, consider
 the statistics $(x_j)_{1\leq j \leq d}=\{ x_{j,k}; k \in \BBz\}_{1\leq j \leq d}$ defined by
\begin{eqnarray}
 x_{j,k} &=& \int_{[0,1]^d} {\bar \phi}_{j,k}^d (t) dX (t) \nonumber \\
& =&  \xi_j \int_{[0,1]} \phi_{k}^1 (t_j) f_j(t_j) dt_j + \epsilon  \eta_{j,k} \nonumber\\
& =& \xi_j \theta_{j,k}  + \epsilon  \eta_{j,k}  \label{transformed_obs},
\end{eqnarray}
where  the  random variables
$\eta_{j,k}=\int_{[0,1]^d} {\bar \phi}_{j,k}^d(t)  dW(t)$  are  i.i.d. real standard Gaussian random variables and  $\theta_{j,k}=\int_{[0,1]} \phi_{k}^1(t_j) f_j(t_j) dt_j$.
Set ${ \theta_j}=(\theta_{j,k})_{k \in \BBz}$ and $\boldsymbol{\theta}= (\theta_j)_{1\leq j \leq d}$.

Thanks  to the  periodic constraints, we may consider $(\phi_k^1)_{k \in \BBz}$ as the standard Fourier basis.
 Then the Sobolev semi-norm of $f_j$ can be expressed in terms of its Fourier coefficients as follows:
 $\|f_j\|_2^{(\tau)} = (\displaystyle{ (2 \pi)^{2 \tau}\sum_{k\in \BBz}} |k|^{2 \tau}\theta_{j,k}^2)^{1/2}$.
Therefore,
 the functional class ${\cal F}_d(\tau, r_\epsilon,b)$ can be equivalently represented as the sequence space $\Theta_d (\tau,r_{\epsilon},b)$:
$$\Theta_d (\tau,r_{\epsilon},b) = \{ \overline{\boldsymbol{\theta}}=( { \theta_1}\xi_1,\ldots,{ \theta_d}\xi_d) : \;  \sum_{j=1}^d \xi_j =d^{1-b};   \;    \; \forall j \in \{1,\ldots,d\},  \theta_j \in \Theta (\tau, r_\epsilon)\},$$
where
$$\Theta (\tau, r_\epsilon) =\{ \theta \in l_2 (\BBz) : \; (2 \pi)^{2 \tau}  \sum_{k\in \BBz}| k|^{2\tau}  \;  \theta_{k}^2 \leq 1;   \sum_{k\in \BBz}
 \; \theta_{k}^2 \geq r_{\epsilon}^2
 \}.$$
The  testing problem  of interest (\ref{PB-Test}) can be rewritten
 in the form
$$
H_0 \; : \; \overline{\boldsymbol{\theta}}  = 0 \;\;\;\;\;\;
 \mbox{ {\rm versus} }  \;\;\;\;\;\;
  H_{1} \; : \;  \overline{\boldsymbol{\theta}} \in \Theta_d(\tau,r_\epsilon,b).$$

Denote by $\P_0$ and    $\P_{\overline{\boldsymbol{\theta}}}$ the distributions  under  the null and alternative hypotheses, respectively.
Also, denote by $\E_0$, $\V_0$,  $\E_{\overline{\boldsymbol{\theta}}}$, and $\V_{\overline{\boldsymbol{\theta}}}$ the expectations and variances  with respect to $\P_0$ and $\P_{\overline{\boldsymbol{\theta}}}$, respectively.
The notation $\P_{{ \theta_j}}$, $\E_{{ \theta_j}}$ and $\V_{{ \theta_j}}$ also will be used:  they are related  to the
distribution of the  observations $x_j=(x_{j,k})_{k \in \BBz}$.

For any test procedure $\psi$, that is, for any function  measurable with respect to the observations and taking its values on the interval $[0,1]$, let
$ \omega(\psi)=\E_0( \psi)$ be the type I error probability and let  $ \beta (\psi,\Theta_d(\tau,r_\epsilon,b))=\displaystyle{\sup_{{\overline{\boldsymbol{\theta}}} \in \Theta_d(\tau,r_\epsilon,b) }}
\E_{\overline{\boldsymbol{\theta}}} (1 - \psi)$ be the  maximal type II error probability over the set  $\Theta_d(\tau,r_\epsilon,b)$.
Also, consider  the total error probability
$ \gamma (\psi,  \Theta_d(\tau,r_\epsilon,b))=\omega(\psi) + \beta (\psi, \Theta_d(\tau,r_\epsilon,b))$,
and denote by  $\gamma$ or $ \gamma ( \Theta_d(\tau,r_\epsilon,b))$ the minimax total error probability
over $\Theta_d(\tau,r_\epsilon,b)$, that is,
\begin{eqnarray} \gamma &= & \gamma (\Theta_d(\tau,r_\epsilon,b))=
 \displaystyle{\inf_{\psi}} \gamma(\psi,\Theta_d(\tau,r_\epsilon,b)), \label{gamma} \end{eqnarray}
where the infimum     is taken over all  test procedures. One can not distinguish between $H_0$ and $H_1$ if
$ \gamma \rightarrow 1$, and
distinguishability occurs if it exists $\psi$ such that
either $\gamma (\psi,  \Theta_d(\tau,r_\epsilon,b)) \rightarrow 0$ or
$\beta (\psi, \Theta_d(\tau,r_\epsilon,b))=o(1)$ once $\psi$
has a given  asymptotic level.

The aim of this paper is to provide separation rates for the alternatives
$\Theta_d(\tau,r_\epsilon,b)$ and to determine statistical
procedures $\psi$ and/or $\psi_\a$ asymptotically of level $\alpha$,
i.e., $\omega(\psi_\a)\le \a+o(1)$, for which these separation
rates are achieved.

By the separation rates we mean
a family $r^\star_\epsilon$ such that \\

$\left\{ \begin{array}{lll}
 \gamma \rightarrow 1  & \mbox{{\rm if }} & \displaystyle{\frac{r_\epsilon}{r_\e^\star}} \rightarrow 0, \\
 && \\

\gamma (\psi,  \Theta_d(\tau,\epsilon,b))\to 0,\quad \mbox{ {\rm and/or} }\quad \forall\ \a\in (0,1) \quad \beta (\psi_\a, \Theta_d(\tau,r_\epsilon,b))  \rightarrow 0  & \mbox{{\rm if }} & \displaystyle{\frac{r_\epsilon}{r_\e^\star}} \rightarrow +\infty.
 \end{array} \right.$\\

By the sharp separation rates, we mean
a family $r_\e^\star$ such that \\

$\left\{ \begin{array}{lll}
  \gamma \rightarrow 1  & \mbox{{\rm if }} & \displaystyle{\limsup \frac{r_\epsilon}{r_\e^\star}} < 1, \\
 && \\

\gamma (\psi,  \Theta_d(\tau,r_\epsilon,b))\to 0,\quad \mbox{ {\rm and/or} }\quad \forall\ \a\in (0,1) \quad \beta (\psi_\a, \Theta_d(\tau,r_\epsilon,b))  \rightarrow 0& \mbox{{\rm if }} & \displaystyle{\liminf \frac{r_\epsilon}{r_\e^\star}} > 1. \end{array} \right.$\\

Typically, asymptotics for models like model (\ref{M1}) are  given as   $\epsilon \rightarrow 0$. However, we are mainly  interested
in high-dimensional settings when $d \rightarrow +\infty$.
Therefore, here and later, asymptotics and  symbols  $o$, $O$, $\sim$ and $\asymp$ are used when $\epsilon \rightarrow 0$ and $d \rightarrow +\infty$,
 except for the cases when
it is explicitly specified, say,  $o_d$ is used when $d \rightarrow +\infty$.
The notation $A\stackrel{\Delta}{=}B$ means that we use notation $A$ for quantity $B$.

\section{Extremal  problem} \label{sec:Extreme}

In this section, we explain what is the quantity $a(r_\e)$ that corresponds to the energy of a signal
in the vectorial case. Only in this section, we assume that the observations have the form
$x_k = \theta_k + \e \eta_k$ for $ k \in \BBz$,  where the $\eta_k$'s are i.i.d. real standard Gaussian random  variables.
The quantity $a(r_\e)$ denotes  the solution  of the extremal problem
\begin{equation} a^2(r_\epsilon) =\frac{1}{2\epsilon^4} \inf_{{\theta} \in
  l_2(\BBz)} \sum_{k \in \BBz} \theta_{k}^4 \mbox{ {\rm subject to } } \left\{ \begin{array}{l}
    (2\pi)^{2 \tau} \sum_{k\in \BBz}|k|^{2 \tau} \theta_{k}^2 \leq 1\\
  \sum_{k\in \BBz} \theta_{k}^2  \geq r_\epsilon^2
\end{array} \right. \label{minimisation}\end{equation}
and characterizes distinguishability in the minimax detection problem
for one-variable functions lying in ${\tilde
S}_\tau$  and separated from the null hypothesis in $L_2$  by
positive values $r_\epsilon$, i.e., for $t\in [0,1]$, $f(t)= \displaystyle{\sum_{k \in
\BBz}} \theta_k \phi_k^1(t)$ with $f \in {\tilde S}_\tau$ and
$\|f\|_2 \geq r_\epsilon$.

Namely, if $a(r_\epsilon)\to 0$ then  the minimax total error probability  $ \gamma(\Theta(\tau, r_\e))
\to 1$,
and if $a(r_\epsilon)\to + \infty$, then $\g(\Theta(\tau, r_\e))\to
0$.

Furthermore, let $\theta^\star\stackrel{\Delta}{=}\theta^\star(r_\e)$ be a sequence in
$l_2(\BBz)$ that provides solution to the extremal  problem
(\ref{minimisation}). Set
\begin{eqnarray}
   w_k(r_\e)&=& \frac 1 2 \frac{(\theta^\star_{k}(r_\e))^2}{a(r_\epsilon)\epsilon^2}, \ k\in\BBz
   \label{weights_ideal}.
  \end{eqnarray}
Suppose that
\begin{equation}\label{assumption1}
a(r_\epsilon)\asymp 1,\quad \sup_{k\in\BBz} w_k(r_\e)=o(1).
\end{equation}
Then,  we  get the sharp asymptotics $$ \g(\Theta(\tau,
r_\e))=2\Phi(-a(r_\epsilon)/2)+o(1). $$

For the reader's convenience, we give a sketch of the proofs of these
results. The proofs are based on the methods and results of Sections 3.1, 3.3, 4.3 in
\cite{IS.02a}. In the vectorial case in hand, we also describe the structure of asymptotically
minimax tests.

In order to obtain lower bounds, we consider the Bayesian hypothesis testing
problem with the product prior distribution on $\theta$, using the symmetric two-point factors:
$\pi=\displaystyle{\prod_{k\in\BBz}\pi_k},\ \pi_k=\frac
12(\delta_{-\theta_k}+\delta_{\theta_k})$ for $\theta\in
\Theta(\tau, r_\e)$, and $\delta$ is the Dirac mass.  Let $\P_\pi$ be the mixture of measures
$\P_\theta$ over $\pi$. Observe that
$$
\frac{d\P_\pi}{d\P_0}((x_k)_{k \in \BBz})=\prod_{k\in\BBz}\frac{d\P_{\pi_k}}{d\P_0}(x_k)=\prod_{k\in\BBz}\exp(-\theta_k^2/2\e^2)\cosh(x_k\theta_k/\e^{2}).
$$
For the sake of simplicity, denote $\displaystyle{\frac{d\P_\pi}{d\P_0} \stackrel{\Delta}{=}\frac{d\P_\pi}{d\P_0}((x_k)_{k \in \BBz})}$. Since $\pi(\Theta(\tau, r_\e))=1$, we have, see Proposition 2.12 in
\cite{IS.02a},
$$
 \gamma(\Theta(\tau, r_\e)) \ge 1-\frac 12 \E_0|d\P_\pi/d\P_0-1|\ge 1-\frac 12
(\E_0(d\P_\pi/d\P_0-1)^{2})^{1/2}=1-\frac 12
((\E_0(d\P_\pi/d\P_0)^2)-1)^{1/2}.
$$
This yields $\gamma(\Theta(\tau, r_\e)) \to 1$ as soon as $\E_0(d\P_\pi/d\P_0)^2\to 1$. Simple
calculations and the inequality $\cosh(x)\le \exp(x^2/2)$ give
\begin{eqnarray*}
\E_0(d\P_\pi/d\P_0)^2=\prod_{k\in\BBz}\E_0
(d\P_{\pi_k}/d\P_0)^2=\prod_{k\in\BBz}\cosh((\theta_k/\e)^2)\le
\exp\l(\frac{1}{2\e^4}\sum_{k\in\BBz}\theta_k^4\r).
\end{eqnarray*}
Therefore, providing the  "asymptotically least favorable prior" of
the type under consideration leads to the problem \nref{minimisation}.

Under assumption \nref{assumption1}, taking the prior based on the
extremal sequence in the problem \nref{minimisation}, one can show
that the Bayesian log-likelihood ratio is asymptotically Gaussian:
$$
\log(d\P_\pi/d\P_0)=\sum_{k\in\BBz}\left(-\frac{(\theta_k^\star)^{2}}{2\e^2}+\log(\cosh(x_k
\theta_k^\star/\e^{2}))\right)=-a^2(r_\e)/2+a(r_\e)\eta_\e+\rho_\e,
$$
where $\eta_\e\to\eta\sim\CN(0,1)$  and $\rho_\e\to 0$ in
$\P_0$-probability. The proof is based on Taylor's expansion, see
Section 4.3.1 of \cite{IS.02a}. This yields the sharp lower bounds.

In order to obtain upper bounds, take a sequence $q=(q_k)_{k\in\BBz}$ such that
$q_k\ge 0,\ \sum_k q_k^2=1/2$, and consider $t_q$, a
centered and normalized (under $\P_0)$ statistic  of a weighted
$\chi^2$-type:
$$
t_q=\sum_{k\in\BBz}q_k\left((\frac{x_k}{\e})^2-1 \right).
$$
Consider also the test procedures  $\psi_{H,q}=\1_{t_q>H}$. Observe that $ \E_0 t_q=0,\
\Var_0 t_q=1,$ and $t_q$ are asymptotically standard Gaussian under
$\P_0$. These observations imply $w(\psi_{H,q})= \Phi(-H)+o(1)$.  Denote by
$\kappa({ \theta},q)$ and $\kappa(q)$ the following
 functions:
\begin{eqnarray}
\kappa({ \theta},q)=\sum_{k \in \BBz} q_k \theta^2_{k},\quad
\kappa(q)= \kappa(\Theta(\tau, r_\e),q)=\inf_{\theta\in \Theta(\tau,
r_\e)}\kappa({ \theta},q). \label{function_kappa}
\end{eqnarray}
Then,
$$
\E_\theta t_q=\e^{-2}\kappa({ \theta},q),\quad \Var_\theta
t_q=1+4\e^{-2}\sum_k q_k^2\theta_k^2=1+O((\max_k q_k)\E_\theta t_q),
$$
and hence, by Chebyshev's inequality,
$\beta(\psi_{H,q},\Theta(\tau, r_\e))\to 0$ when
$\e^{-2}\kappa(q)\to +\infty$ and $H\le c\e^{-2}\kappa(q),\ c\in
(0,1)$. Under assumption \nref{assumption1}, one can check that the statistic $\hat
t_q=t_q-\E_\theta t_q$ is asymptotically standard Gaussian under
$\P_\theta$ such that $\E_\theta t_q=O(1)$. Therefore
$$
\beta(\psi_{H,q},\Theta(\tau, r_\e))\le
\Phi(H-\e^{-2}\kappa(q))+o(1).
$$
In order   to determine ``asymptotically the best sequence'' $(q_k)_{k \in \BBz}$, it suffices to find  a
solution of the following maximin problem:
\begin{equation}\label{maxmin}
\tilde a(r_\e)=\e^{-2}\sup_{\sum_{k}q_k^2=1/2, q_k\ge 0}\kappa(q).
\end{equation}
First, we change the variables for $v=(v_k)_{k\in\BBz}$ and $(p_k)_{k\in\BBz}$, where
$v_k=\theta_k^2/\sqrt{2}$, $p_k=\sqrt{2}q_k$. Then, by convexity of
the set
\begin{equation}\label{V+}
V^+ 
=\{v \in l_1(\BBz): v_k\ge 0;\; (2 \pi)^{2 \tau} \sum_{k\in
\BBz}k^{2 \tau} v_k \leq 2^{-1/2}; \;
  \sum_{k\in \BBz} v_k  \geq  2^{-1/2} r_\epsilon^2 \},
\end{equation}
and using the minimax theorem, we get
\begin{eqnarray*}
\tilde a(r_\e) &=&\e^{-2}\sup_{\sum_k p_k^2=1, p_k\ge 0}\inf_{v\in
V^+}\sum_kp_kv_k =\e^{-2}\sup_{\sum_k p_k^2\le
1, p_k\ge 0}\inf_{v\in V^+}\sum_kp_kv_k\\
& =&\e^{-2}\inf_{v\in V^+}\sup_{\sum_k p_k^2\le 1, p_k\ge
0}\sum_kp_kv_k =\e^{-2}\inf_{v\in V^+}(\sum_k
v_k^2)^{1/2}\\
&=&\frac{1}{\sqrt{2}\;\e^2}\inf_{\theta\in\Theta(\tau,
r_\e)}(\sum_k\theta_k^4)^{1/2}=a(r_\e).
\end{eqnarray*}
Thus, asymptotically the best sequence $(q_k)_{k \in \BBz}$ is the sequence $w(r_\e)\stackrel{\Delta}{=}( w_k(r_\e))_{k \in \BBz}$ of the form
\nref{weights_ideal}, and the value of the problem \nref{maxmin}
coincides with the value of the problem \nref{minimisation}. Setting
$H=a(r_\e)/2$, we get the upper bounds and the structure of
asymptotically minimax tests.

Note that the above evaluations  entail  (see also Proposition 4.1 in
\cite{IS.02a}) that
\begin{eqnarray}
\inf_{{ \theta} \in
  \Theta (\tau,r_{\epsilon})} \frac{1}{\epsilon^2}\kappa(\theta,w(r_\e))  \geq  a(r_\epsilon).
  \label{PB_minim} \end{eqnarray}

Moreover if $\displaystyle{(\sum_{k \in \BBz}\theta_k^{2})^{1/2}}$ is larger than $r_\e$, then
$\kappa(\theta,w(r_\e))$ becomes rather large. Namely, let us denote
$$
\kappa(r_\e,B)=\inf_{\theta\in\Theta(\tau,Br_\e)}\kappa(\theta,w(r_\e)), \; B >0
$$
\begin{proposition}\label{P1}
Let $B\ge 1$, then
\begin{equation*}\label{ext1}
\frac{1}{\epsilon^2}\kappa(r_\e,B)\ge B^2a(r_\epsilon).
\end{equation*}
\end{proposition}

\noindent
{\it Proof of Proposition \ref{P1}.} \\ Set $ \Theta(\tau,A,r_\e)=\{\theta\in
l_2(\BBz): (2\pi)^{2\tau}\sum_{k\in\BBz}|k|^{2\tau}\theta_k^2\le
A^2, \;\ \sum_{k\in\BBz}\theta_k^2\ge r_\e^2\}$, $A>0$.  Since
$\Theta(\tau,Br_\e)\subset\Theta(\tau,B, Br_\e)$, we have
\begin{eqnarray*}
\inf_{\theta\in\Theta(\tau,Br_\e)}\kappa(\theta,w(r_\e))\ge
\inf_{\theta\in\Theta(\tau,B, Br_\e)}\kappa(\theta,w(r_\e))=
B^2\inf_{\theta\in\Theta(\tau,r_\e)}\kappa(\theta,w(r_\e))\geq B^2\e^2a(r_\epsilon),
\end{eqnarray*}
where the last inequality follows from
 (\ref{PB_minim}). This completes the proof.

\medskip

The solution of the extremal problem (\ref{minimisation}) is
obtained in Ingster and Suslina \cite{IS.02a}, Section 4.3.
Adapting the derivations on pages 146--147 of  Section 4.3.2. in \cite{IS.02a} to our case,
we set $c_3=\displaystyle{\frac{1}{4 \tau}} \;   B(a,b) $, $c_2 = \displaystyle{\frac{1}{4 \tau}} \;  B(b,c)$ and $c_0=
\displaystyle{\frac{1}{8\tau}} \;   B(a,d)$, where $B(\cdot,\cdot)$ is the Euler Beta function, $a=
\displaystyle{\frac{1}{2\tau}}$,
$b=1+\displaystyle{\frac{1}{2\tau}}$, $c=2$ and $d=3$.
\begin{lemma} \label{sol_min}
The solution of the extremal problem (\ref{minimisation}) is given by
\begin{eqnarray}
a(r_\epsilon)  \sim (c_1(\tau))^{1/2}  \; r_\epsilon^{2 +1/(2\tau)} \; \epsilon^{-2} \; \mbox{ {\rm as } } \; r_\epsilon \rightarrow 0 \label{sol_a}, 
\end{eqnarray}
\begin{eqnarray}
 \mbox{ {\rm where } }  \;\;  c_1(\tau) & =&  c_0  \pi c_2^{-2} (\frac{c_2}{c_3})^{(4 \tau +1)/2 \tau} \mbox{ {\rm is a positive constant.} }  \label{constant-c1}
\end{eqnarray}
\end{lemma}
\noindent

\begin{remark} \label{rm-reps-0}
One must note that  $r_\epsilon \rightarrow  0$ is the only condition we need to obtain the asymptotic solution of (\ref{sol_min}).
In particular, it is not required that  $\epsilon \rightarrow 0$  and Lemma \ref{sol_min} is  valid whatever
the  value of  $\epsilon>0$ is.
\end{remark}

\noindent {\it Sketch of proof of Lemma \ref{sol_min}.} \\
Following Chapter 4  in \cite{IS.02a}, observe that by setting $v_k = \theta_k^2/{\sqrt 2 }$ for all $k \in \BBz$, one can transform the minimization problem under constraints (\ref{minimisation}) into the following one:
\begin{eqnarray*} v^2_\epsilon& =& \inf_{{(v_k)_{k \in \BBz} \in V^+}} \sum_{k \in \BBz} v_k^2, \end{eqnarray*}
  where $V^+$ is defined by  equation (\ref{V+}).
The space $l_1^+(\BBz)$ contains  non-negative  sequences lying in $l_1(\BBz)$.
Note that  $v_\epsilon^2= \epsilon^4 a^2 (r_\epsilon)$.  The convexity of the set
$ V^+$ assures the uniqueness of
$v_\epsilon^2$.
In order to determine the solution, rewrite as in Section 4.3. in \cite{IS.02a} the sequence $(v_k)_{k \in \BBz}$ as follows:
$v_k=v_0 \zeta(k/m)$, where   $\zeta(y)=   (1-|y|^{2 \tau} ) \I_{(|y| \leq 1)}$ and $m>0$. 
By using the Lagrange multipliers rule, it is possible to obtain
the following relations, as $r_\epsilon \rightarrow 0$ and $m \rightarrow + \infty$:
\begin{equation} \label{Lagrange}
c_{3} v_0 m \;\sim \; 2^{-1/2} r_\epsilon^2,   \; \;
v_\epsilon^2  \; \sim  \; c_0  v_0^2 m, \; \;
c_2 v_0 m^{2 \tau +1} \; \sim \; 2^{-1/2}  (2 \pi)^{-2 \tau},
\end{equation}
which entail the existence of
$v^2_\epsilon$ satisfying $v^2_\epsilon \sim c_1(\tau) r_\epsilon^{4+1/\tau}$,
   and thus $a^2(r_\epsilon) \sim c_1(\tau) \epsilon^{-4} r_\epsilon^{4+1/\tau}$.

If  $r_\epsilon \rightarrow 0$, then the first and  second relations in (\ref{Lagrange}) entail that
\begin{eqnarray}
v_0 \asymp v_\epsilon^2 r_\epsilon^{-2} \asymp r_\epsilon^{2+1/\tau}, \label{sup-vk}
\end{eqnarray}
which implies that $m \rightarrow +\infty$ since the third relation in (\ref{Lagrange}) yields $m\asymp v_0^{-1/(2\tau+1)} \asymp r_\epsilon^{-1/\tau}$.
\begin{remark} \label{rmv0} The form of
function $\zeta$ and relation (\ref{sup-vk}) imply that  ${\displaystyle{\sup_{k }}}\, v_k \leq v_0=o(1)$. 
\end{remark}

\section{Main results} \label{sec:Main}
Depending on the values of
$b$, we distinguish between two types of sparsity: the moderate sparsity case with $b\in(0,1/2]$ and the high
sparsity case with $b\in(1/2, 1)$.  In each case, although being of different types, the ``best''
test procedures that achieve the separation rates are based on the $\chi^2$-type statistics $(t_j)_{1 \leq  j \leq  d}$
determined in the same way as the ``best statistic'' $t_q$ of a weighted $\chi^{2}$-type in Section \ref{sec:Extreme}.

Let us introduce a general version of the $\chi^2$-type statistics of interest. For $j$ in $\{1, \ldots, d\}$, put
\begin{eqnarray}
t_j &=& \sum_{k\in \BBz} w_k \left(
(\frac{x_{j,k}}{\epsilon})^2-1 \right) \label{tj},
\end{eqnarray}
 where $(w_k)_{k \in \BBz}$ is the sequence of
  weights such that $w_k \geq 0$ for all $k$ in $\BBz$ and  $\sum_{k \in \BBz} w_k^{2}=\frac{1}{2}$.
    Set also
  \begin{eqnarray}
  t_{j,k} & =&  w_k \left( (\frac{x_{j,k}}{\epsilon})^2-1 \right)\label{tjk},
  \end{eqnarray}
so that $t_j  = \displaystyle{\sum_{k\in \BBz}t_{j,k}}$.



Recall that $T_d=\sqrt{ \log d}$ (see Section \ref{sec:gaussian-vector}).  Similarly to (\ref{HC-DJ04}) and  for
any $u \in (0, \sqrt{2}]$, let us  define  the statistics
$L(u)$ on which the Higher-Criticism type test procedure is  built:
\begin{eqnarray}
L(u)= C_u \sum_{j=1}^d (\I_{(t_j > u T_d)} - {\tilde \Phi}_0(
u T_d)), \label{Inter-Hig}
\end{eqnarray}
where 
   \begin{eqnarray}{\tilde \Phi}_0(x)&=&\P_0 (t_{j} > x) \label{Phi} \\
   C_u &=& (d {\tilde \Phi}_{0}( uT_d)
(1-{\tilde \Phi}_{0}( uT_d)))^{-1/2}. \label{Cu}
\end{eqnarray}

Taking into account the sparsity condition, we consider a particular sequence of weights $(w_k(r^\star_\e))_{k \in \BBz}$  defined by equation \nref{weights_ideal} 
with  $r^\star_\e\stackrel{\Delta}{=}r^\star_\e(b)$ being the separations rates depending on
 $b$  in $(0,1)$. Then, for all $j \in \{1,\ldots,d\}$,
we consider the statistics $t_{j,b}$ as in (\ref{tj}) with the weight sequence $(w_k(r^\star_\e))_{k \in \BBz}$, that is,
$$
t_{j,b} =\sum_{k\in \BBz} w_k(r^\star_\e) \left(
(\frac{x_{j,k}}{\epsilon})^2-1 \right) .
$$
Also, denote by  $t_b$  the normalized empirical
mean of the $t_{j,b}$'s:
\begin{eqnarray}
t_b= \frac{1}{{\sqrt d}} \; \sum_{j=1}^d t_{j,b}.  \label{stat_moderate}
\end{eqnarray}
Similarly, replacing $t_j$ by $t_{j,b}$, consider the statistics $L(u,b)$, $C_{u,b}$, and ${\tilde \Phi}_{0,b}$  defined by equations (\ref{Inter-Hig}), (\ref{Cu}) and
(\ref{Phi}) respectively, that is,
\begin{eqnarray}
L(u,b)&= &C_{u,b} \sum_{j=1}^d (\I_{(t_{j,b} > u T_d)} - {\tilde \Phi}_{0,b}(
u T_d)),  \label{Lub} \\
 C_{u,b}&=&(d {\tilde \Phi}_{0,b}( uT_d)
(1-{\tilde \Phi}_{0,b}( uT_d)))^{-1/2},\nonumber\\
{\tilde \Phi}_{0,b}(x)&=&\P_0 (t_{j,b} > x).\nonumber
 \end{eqnarray}



\subsection{Moderate sparsity}
In case of moderate sparsity, for any $\alpha \in (0,1)$,
consider the  $\chi^2$-type test procedure:
\begin{eqnarray}
\psi^{\chi^2}_\alpha \stackrel{\Delta}{=}\psi^{\chi^2}_{\alpha,b} & =& \I_{( t_b >
T_\alpha)}, \label{Test_Moyenne_Emp}
\end{eqnarray}
where $t_b$ is defined in  (\ref{stat_moderate}) and
$T_\alpha$ is 
 the  $(1-\alpha)$-quantile  of a real standard Gaussian random variable.

\begin{theorem} \label{Moderate_case} 
Assume that $r_\epsilon \rightarrow 0$ and let  $a(r_\epsilon)$ be given by
 (\ref{sol_a}). Then, 
  the following results hold true. 
\begin{itemize}
\item  (i) Lower bound. \begin{description}
\item If $a(r_\epsilon) d^{1/2 -b} =o(1)$, then
 $\gamma  \rightarrow 1$.
\item
If  $a(r_\epsilon) d^{1/2 -b}=O(1)$,   then $\lim \inf \gamma >0$.
\end{description}
\item (ii) Upper bound. Let $r_\e^\star=r_\e^\star(b)$ be determined by the relation $a(r_\epsilon^\star)\asymp
 d^{b-1/2}$ and $\psi^{\chi^2}_\alpha$ be defined by
 (\ref{Test_Moyenne_Emp}).
 Then,
\begin{description}
\item Type I error: $\forall \alpha \in (0,1)$, $\omega ( \psi^{\chi^2}_\alpha) =\alpha + o(1)$.
\item Type II error: if  $a(r_\epsilon) d^{1/2 -b} \rightarrow  +\infty$,  then $\beta(\psi^{\chi^2}_\alpha,\Theta_d(\tau,r_\epsilon,b))=o(1)$.
\end{description}

\end{itemize}
\end{theorem}
 \begin{remark}
 Note that we obtain the same detection boundaries  as in the vectorial case (see Section \ref{sec:Bgk-Seq-Sp}):  the areas of distinguishability and non-distinguishability depend on
the limit of $d^{1/2 -b} a(r_\epsilon)$.
 The  condition   $d^{1/2 -b} a(r_\epsilon) \asymp a(r_\epsilon)/a(r_\epsilon^{\star}) \rightarrow +\infty$
   is equivalent to
 $  r_\epsilon/ r_\epsilon^{\star} \rightarrow +\infty$ where by (\ref{sol_a})
 \begin{eqnarray} r_\epsilon^{\star}\asymp (\epsilon^4 d^{2b -1})^{\tau/(4 \tau +1)}. \label{sep-rates-M} \end{eqnarray}
In order to use Lemma \ref{sol_min},  the condition $r_\epsilon \rightarrow 0$ is required.
Note that the requirement $r_\epsilon^\star \rightarrow 0$ is always fulfilled for $b \in (0,1/2)$ whatever the value of $\epsilon>0$ is 
as soon as $d \rightarrow +\infty$. For  $b=1/2$, the condition $r_\epsilon^\star \rightarrow 0$ holds when
$\epsilon \rightarrow 0$.
\end{remark}

\subsection{High sparsity}
Let us define  the Higher-Criticism type test procedure. Let
$r^\star_\e=r^\star_\e(b)$ be determined by the relation $a(r^\star_\e)\sim
\varphi(b)T_d$, where $\varphi(b)$ is given by \nref{phi}. Set
$u(b)=\min(2\varphi(b),\sqrt{2})$, i.e., $u(b)=2\varphi(b)$ for
$b\in (1/2,3/4]$, and $u(b)=\sqrt{2}$ for $b\in (3/4,1]$. Consider the test
\begin{eqnarray*}
\psi^{L} & = & \I_{\{ \displaystyle{\max_{1 \leq l \leq N-1}}
L(u_l,b_l)>H\}},\quad u_l=u(b_l) \label{Test_Max1},
\end{eqnarray*}
where  the function  $L$ is defined in (\ref{Lub}) and
$(b_l)_{1 \leq l \leq N}$ consists  of a regular grid  on  $(1/2, 1]$,
that is, $b_l=1/2 + l \delta$, where $\delta$ is a positive parameter that
satisfies $\delta=o_d(1)$, $T_d \delta \rightarrow + \infty$ and $N
\delta = 1/2$. This entails that $N=O_d(\delta^{-1})$ and thus
$N=o_d(T_d)$.   Take a positive $H$ such that $H \sim (\log d)^C$ for
some
positive constant $C$ satisfying  $C > \frac 1 4$. 

For a constant $D>\sqrt{2}$, consider also the test
\begin{eqnarray*}
\psi^{max} & = & \I_{\{ \displaystyle{\max_{1 \le j\le d}\max_{1\le
l\le N} t_{j,b_l}>D T_d}\}}
\label{Test_Max2}.
\end{eqnarray*}
Finally, combining $\psi^{L}$ and  $\psi^{max}$, we define  the test procedure
\begin{eqnarray}
\psi^{HC}=\psi^{L}\psi^{max},
\label{Test_Max}
\end{eqnarray}
that rejects $H_0$ if both $\psi^{L}$ and
$\psi^{max}$ reject $H_0$.

For the high sparsity case, not only separation rates but also sharp
asymptotics are obtained; two ranges of $b$ should be distinguished:
the range of $b$ in $(1/2, 3/4]$, called the intermediate
sparsity case, and  the range of $b$ in $(3/4,1)$, called the
highest sparsity case.

\begin{theorem} \label{Intermediate_case}
Assume that  $r_\epsilon \rightarrow 0$ and that
$\log d = o (\epsilon^{-2/(2 \tau +1)})$.  Let $a(r_\epsilon)$ be given by (\ref{sol_a})
and let $\varphi$ be given by  (\ref{phi}).

\begin{itemize}
\item (i) Lower bound. If
$\limsup a(r_\epsilon)/T_d < \varphi(b) $, then  $\lim \inf  \gamma  \rightarrow 1$.
\item (ii)  Upper bound: errors  of $\psi^{HC}$ defined by (\ref{Test_Max}). 
\begin{itemize}
\item Type I error:
$\omega (\psi^{HC})=o(1)$.
\item Type II error:
if $\liminf a(r_\epsilon)/T_d > \varphi(b)$, then
$ \beta(\psi^{HC},\Theta_d(\tau,r_\epsilon,b))=o(1)$.
\end{itemize}
\end{itemize}
\end{theorem}


%
%
%

\begin{remark} \label{rm:appro}
\begin{itemize}
\item Set  $ a (r_\epsilon^{\star})=T_d \varphi (b)$.
In our sparse functional framework, the distinguishability conditions are the same as for a
$d$-dimensional sparse vector (see, e.g., \cite{IS.02b}), with the only difference that in our case the assumption $\log d =o(\epsilon^{-2/(2 \tau +1)})$ is required. Under this assumption, the result of Theorem  \ref{Intermediate_case}
 means that   distinguishability is impossible if
$ \limsup a(r_{\epsilon})/ a (r_\epsilon^{\star})< 1$ and it is possible if  $ \liminf  a(r_{\epsilon})/ a (r_\epsilon^{\star})> 1$.  Due to (\ref{sol_a}), these conditions provide sharp separation rates since they are equivalent to  $ \limsup r_{\epsilon}/ r^\star_{\epsilon}< 1$ and $\liminf r_{\epsilon}/ r^\star_{\epsilon}> 1$, respectively,  where
\begin{eqnarray}
r_\epsilon^\star&\sim &
( \epsilon^4  \; T_d^2 \; (c_1(\tau))^{-1} \varphi^2(b))^{\tau/(4 \tau +1)},  \label{sep-rates-H}
\end{eqnarray}
and  $c_1(\tau)$ is  defined by (\ref{constant-c1}).  Note that the condition $r_\epsilon^\star \rightarrow 0$ is fulfilled
under the assumption $\log d =o(\epsilon^{-2/(2 \tau +1)})$.
\\
The values $r_\epsilon^\star$ mark the border between the areas of  distinguishability and   non-distinguishability. Indeed,
for $r_\epsilon \rightarrow 0$  such that
$\lim \sup r_\epsilon/r_\epsilon^{\star} < 1$,  the alternatives  separated from the null hypothesis by $r_\epsilon$ are not distinguishable  and,  on the other side,
for  $r_\epsilon \rightarrow 0$  such that   $\lim \inf r_\epsilon/ r_\epsilon^{\star} > 1$, the alternatives separated from the \textbf {null hypothesis} by $r_\epsilon$ are
distinguishable.
\item  Actually, the assumption $\log d =o(\epsilon^{-2/(2 \tau +1)})$ is equivalent  to \begin{eqnarray}
(r_{\epsilon}^{\star})^{1/(2\tau)} T_d =o(1), \label{necc-Lag-sup}                                                                                            \end{eqnarray}
which  is required   when dealing with the asymptotic behavior of the tail distribution of $t_{j,b}$ (see Lemma \ref{resu_tech})  since
$T_d \; \displaystyle{\sup_k} w_k(r_\e^\star)\leq (r_{\epsilon}^{\star})^{1/(2\tau)} T_d$.  Relation
(\ref{necc-Lag-sup})  follows from the relations in (\ref{Lagrange}).  Concerning the lower bound, condition (\ref{necc-Lag-sup})  is necessary
when we  evaluate  the second moment of the Bayesian likelihood ratio under the null hypothesis.
\item Note that the condition $\log d =o(\epsilon^{-2/(2 \tau +1)})$
is essential for $b\in (1/2,1)$. Namely, it follows from Theorem 2 in
\cite{IL.03} that if $\lim\inf \; (\log d \; \epsilon^{2/(2 \tau +1)})>0$,
then the separation rates are of the form $r^\star_\e=\e\sqrt{\log d }$
for any $b\in (1/2,1)$. Observe that if $\log d \ge c\e^{-2}$ for
some $c>0$, then the separation rates are bounded away from zero,
i.e., it is impossible to detect  functions lying in
$\Theta_d(\tau,r_\epsilon,b)$ with small enough $r_\e>0$.
\end{itemize}
\end{remark}

\begin{remark} {\sc Adaptation.}\label{adapt}

In the high sparsity case, a family of test procedures $\psi^{HC}$ provides the distinguishability for all $b\in
(1/2,1)$. Moreover,  it follows from the proofs that our result is uniform
over $b\in (1/2+\rho, 1-\rho)$ for any $\rho\in (0,1/4)$, i.e., the results are  adaptive
over $b\in (1/2+\rho, 1-\rho)$ for any $\rho\in (0,1/4)$, without a loss in separation rates.

For the moderate sparsity case, it is worth noting that
the family of test procedures $\psi_\alpha^{\chi^2}=\psi_{\alpha,b}^{\chi^2}$
depends on $b\in (0,1/2]$ since the sequence of weights
$w(r^\star_\e(b))$ does.  It is shown in Theorem 3 of \cite{IL.03}
that ``adaptive'' separation rates for
unknown $b\in (0,1/2)$ are of the form $r_\epsilon^{\star}\asymp
(\epsilon^4 d^{2b -1} \;\log\log d )^{\tau/(4 \tau +1)}$, i.e., the adaptive case leads to an unavoidable   $\log\log$-loss   in separation rates  compared to non-adaptive setting.
Using the Bonferroni method, it is   possible to prove that  the  test procedures based on
a grid of tests of the form $\psi_{\alpha_d, b_l}^{\chi^2}$
are adaptive rate-optimal test procedures.
Since  this result is  similar to the one stated in \cite{IL.03},  we omit it.

\end{remark}

\section{Extended problem} \label{sec:Extended}
In this section, we generalize  the   hypothesis testing problem stated in (\ref{PB-Test})  to more general alternatives.
The null hypothesis $H_0$ is still characterized by
some constant $const_0$ and,  as in (\ref{PB-Test}), under the alternative, the signal function $f$ is, up to some constant, equal to $f^1$, i.e., $f=const_1 + f^1$.
The additive sparse structure on $f^1$ is still assumed, i.e., $f^1 \in {\cal F}_{d,b}$, as well as  every component $f^1_j$ is assumed $1$-periodic and
 orthogonal to a constant (recall that for any $t\in \left[0,1\right] ^{d}$
$f^1(t)= \sum^{d}_{j=1} \xi_j f^1_j(t_j)$ where $\xi_j \in \{0,1\}$ and $t_j \in \left[0,1\right]$  for any $j \in \{1, \ldots,d\}$ ). We then denote by ${\tilde {\cal F}}_{d,b}$ the set of signal functions in ${\cal F}_{d,b}$
 whose components  are $1$-periodic and orthogonal to a constant.
Rather than imposing smoothness constrains component-wise,
we now study the alternative classes for which the smoothness and separation conditions  are expressed in terms of the whole signal function $f^1$.
In other words, the main difference between
the extended and initial detection problems
is that the   distinguishability problem is studied with respect to a global signal.

Then, given the  alternatives that include signal functions $f$ as in (\ref{PB-Test}), where $f^1$ belongs to the functional class
 ${\cal  F}_d^{ext}(\tau,L,r_{\epsilon},b)$, the testing problem of interest is stated as follows:
\begin{eqnarray}H_0\; :\; f= const_0 \;\;\;\;\; \mbox{{\rm
versus }} \;\;\;\;\;  H_{1} \;: \; f=const_1+f^{1}, \; f^{1} \in {\cal
F}^{ext}_d(\tau,L, r_\epsilon, b),\label{PB-Test-Ex}\end{eqnarray}
where
\begin{eqnarray*}
{\cal  F}_d^{ext}(\tau,L,r_{\epsilon},b) & =& \left\{  f^1 \in {\tilde {\cal F}}_{d,b} :\|f^1\|_2\geq r_{\epsilon},\; \|f^1\|_2^{(\tau)} \leq L \right\},
\end{eqnarray*}
in which
$(\|f^1\|_2^{(\tau)})^2=\sum_{j=1}^d \xi_j (\|f^1_j\|_2^{(\tau)})^2$. Due to the  periodic constraint,  we consider the standard Fourier basis.
 This allows to express the semi-norm $\displaystyle{\|\cdot \|_2^{(\tau)}}$  in terms of Fourier coefficients.  As
in Section \ref{sec:Trans}, we then transform
  the functional space ${\cal F}_d^{ext}(\tau,r_{\epsilon},L,b)$ to the sequence space  $\Theta_d^{ext}(\tau,L,r_\epsilon,b)$, which consists of sequences $\overline{{\boldsymbol \theta}}=(\xi_j  \theta_{j,k})_{j,k}$ such that
\begin{eqnarray*}
\sum_{j=1}^d \xi_j = d^{1-b}=K, \\
\;  \sum_{j=1}^d \xi_j  (2 \pi)^{2\tau} \sum_{k\in \BBz}| k|^{2\tau}  \; \theta_{j,k}^2 \leq L^2 , \\
\; \; \sum_{j=1}^d   \xi_j \sum_{k\in \BBz}
 \; \theta_{j,k}^2 \geq  r_{\epsilon}^2.
\end{eqnarray*}

Note that if  $L^2= K$ and ${\tilde r}_{\epsilon}^2=K r_\epsilon^2$,  then we have
\begin{eqnarray*}
 \Theta_d^{ext}(\tau,L,{\tilde r}_\epsilon,b) \supset  { \Theta}_d(\tau,r_\epsilon,b).\label{link-ext}
\end{eqnarray*}
This implies that the results on the lower bound continue to hold for  $\Theta_d^{ext}(\tau,L,{\tilde r}_\epsilon,b)$ with the separation rates    $({\tilde r}_\epsilon^\star)^2 =
 K (r_\epsilon^\star)^2$, where $r_\epsilon^\star$ is defined by either  (\ref{sep-rates-M}) or (\ref{sep-rates-H}) depending on  the values  of $b$.
Here, the quantity of interest is  ${\tilde a}(r_\epsilon)$, the solution of the following extremal problem:
  \begin{equation} {\tilde a}^2(r_\epsilon) =\frac{1}{2\epsilon^4} \inf_{\overline{{\boldsymbol \theta}} \in
  l_2} \displaystyle{\sum_{j=1}^d} \xi_j  \;\sum_{k \in \BBz} \theta_{j,k}^4 \mbox{ {\rm subject to } } \left\{ \begin{array}{l}
     \displaystyle{\sum_{j=1}^d} \xi_j =d^{1-b}=K \\
    \displaystyle{\sum_{j=1}^d}  \xi_j(2 \pi)^{2\tau} \sum_{k\in \BBz}|k|^{2 \tau} \theta_{j,k}^2 \leq K\\
 \displaystyle{\sum_{j=1}^d} \xi_j \sum_{k\in \BBz} \theta_{j,k}^2  \geq K r_\epsilon^2
\end{array} \right. \label{minimisation_Ext}\end{equation}

As follows from  Section 4.3 in  \cite{IS.02a},  the solution of the extremal problem (\ref{minimisation_Ext}) is given by 
\begin{eqnarray*}
{\tilde a}(r_\epsilon) & \sim & (c_1(\tau))^{1/2} K r_\epsilon^{2+1/(2 \tau)} \epsilon^{-2}\; \mbox{ {\rm as } } \; r_\epsilon \rightarrow 0, \label{sol-min-ext}
\end{eqnarray*}
where $c_1(\tau)$ is defined in (\ref{constant-c1}). That is, ${\tilde a}(r_\epsilon) =Ka(r_\epsilon)$,
where $a(r_\epsilon)$ is the solution (\ref{sol_a}) of the extremal problem (\ref{minimisation}).

\begin{remark}
Consider the  function $\kappa$ defined by  (\ref{function_kappa}),  for which  the sequence of weights $w(r_\e)=(w_k(r_\e))_{k}$ is defined as in  (\ref{weights_ideal}).
 Then we obtain from (\ref{PB_minim}) that
\begin{eqnarray}
\inf_{\overline{{\boldsymbol \theta}} \in
  \Theta_d^{ext}(\tau,K^{1/2}, K^{1/2}r_\epsilon,b)} \frac{1}{\epsilon^2}  \sum_{j=1}^d \xi_j
  \kappa (\theta_j,w(r_\e))& \geq &{\tilde   a}(r_\epsilon)= K
  a(r_\epsilon),
 \label{kappa-ext}
\end{eqnarray}
and similarly to Proposition \ref{P1} for any $D\ge 1$,
\begin{eqnarray}
\inf_{\overline{{\boldsymbol \theta}} \in
  \Theta_d^{ext}(\tau,K^{1/2}, D K^{1/2}r_\epsilon,b)} \frac{1}{\epsilon^2}  \sum_{j=1}^d \xi_j
  \kappa (\theta_j,w(r_\e))& \geq & D^2{\tilde   a}(r_\epsilon)= D^2K
  a(r_\epsilon).
 \label{kappa-ext_1}
\end{eqnarray}
\end{remark}

Now, as in Section \ref{sec:Trans}, with the use of the orthonormal system,  instead of considering the random process $X(t)$ defined in  model (\ref{M1}),
we observe  a family of  random sequences
$ (x_{j,k})_{ k \in \BBz,  j \in \{1,\ldots,d\}}$  defined by (\ref{transformed_obs}).
Finally,  the remained question is: do the families of test  procedures $ \psi^{\chi^2}_\alpha$ given by (\ref{Test_Moyenne_Emp}) and $\psi^{HC}$ given by  (\ref{Test_Max})  provide distinguishability?
The answer is affirmative and is given below.
%
Note that it is then sufficient to study the type II error probability of these tests since their type I error probability has been already studied for the
hypothesis testing problem (\ref{PB-Test}).

\begin{theorem} \label{thm:extended-problem}
Assume that $r_\epsilon \rightarrow 0$ and let $ a(r_\epsilon)$ and $\varphi$ be given by (\ref{sol_a}) and (\ref{phi}), respectively.
 Then,  the following results hold true. 
\begin{itemize}
\item (i) {\sc Moderate sparsity}-Type II error probability of $\psi^{\chi^2}_\alpha$ defined by (\ref{Test_Moyenne_Emp}).

 If $a(r_\epsilon) d^{1/2 -b} \rightarrow +\infty$, then
$\beta( \psi^{\chi^2}_\alpha,\Theta^{ext}_d(\tau,K^{1/2}, K^{1/2}r_\epsilon,b))=o(1)$.
\item (ii) {\sc  High sparsity}-Type II error probability of    $\psi^{HC}$ defined by (\ref{Test_Max}).

 Assume  that $\log d = o (\epsilon^{-2/(2 \tau +1)})$. \\
If $ \liminf a(r_\epsilon)/T_d  > \varphi (b) $, then
$ \beta(\psi^{HC}, \Theta^{ext}_d(\tau,K^{1/2}, K^{1/2}r_\epsilon,b)) =o(1)$.
%
\end{itemize}
\end{theorem}

\begin{remark}
One should note that the detection boundaries are the same for the hypothesis testing problems
(\ref{PB-Test}) and (\ref{PB-Test-Ex}), the initial one and
 its generalization.
\end{remark}
%

\section{Proofs} \label{sec:Proofs}
Proofs  of our main results  require some preliminary results that are stated below both under the null and alternative hypotheses.
Specifically,  we establish  asymptotic tail distributions of the test statistics in hand and find their first and second moments.

\subsection{Properties of test statistics}

In this section, we consider the statistics $t_j$ defined by
(\ref{tj}) with any  sequence of weights $w=(w_k)_{k\in\BBz}$ such that $w_k\ge 0, \forall k\in \BBz$ and $\sum_k w_k^2=1/2$. Therefore the quantities  $L(u)$, $C(u)$, and  ${\tilde \Phi}_0$  are those defined by
(\ref{Inter-Hig}), (\ref{Cu}) and (\ref{Phi}).

\begin{proposition} \label{GD} {\sc Asymptotic tail distribution  of $t_j$ defined by
(\ref{tj}).} $\;$

   Assume $T \max_k w_k =o(1)$, then
\begin{eqnarray*}
\log \P_0( t_j > T)  &\sim &   -\frac{T^2}{2}  \;    \mbox{ {\rm  as $T \rightarrow + \infty$}},  \label{Grande_Dev_0} \\
  & & \\
\log \P_{{ \theta_j}} (t_{j} > T)  &\sim  &-\frac{(T-\E_{\theta_j} (t_j))^2}{2 },  \mbox{ {\rm  as $(T-\E_{\theta_j} (t_j))  \stackrel{T \rightarrow +\infty}{\longrightarrow} + \infty$}}. \label{Grande_Dev_1}
\end{eqnarray*}

\end{proposition}

\vspace{0.3cm}
\noindent {\it Proof of Proposition \ref{GD}.} \\
%
We consider only the distribution $\P_{\theta_j}$ since  $\P_0$ is a particular  case of $\P_{\theta_j}$. The proof
consists of bounding  $\P_{\theta_j} (t_j >T)$ from above and below. This is done by using the cumulant-generating function of $t_j$ under  $\P_{\theta_j}$ which is defined by  $\phi_{\theta_j}(h)= \log( \E_{\theta_j}(\exp(h  t_{j})))$ for any  $h$. Let us consider only positive $h$ and let us introduce
a new family of probability  measures $\P_{\theta_j,h}$ such that $\displaystyle{\frac{d \P_{\theta_j,h}}{d \P_0}}=\exp( h  t_j) \;\exp(-\phi_{\theta_j}(h))$.
This yields
\begin{eqnarray}
 \P_{\theta_j} (  t_j > T)& =& \E_{\theta_j,h} [ \I_{(t_j > T)} \exp(-(h t_j- \phi_{\theta_j} (h)))] \nonumber\\
& =& \exp(-(hT -\phi_{\theta_j}(h) )) \; \E_{\theta_j,h} [ \I_{(t_j > T)} \exp(-h( t_j- T))]. \label{equa}
 \end{eqnarray}
Let us start with the upper bound.  \\

\noindent
{\it Upper bound.}
The second term o the right-hand side of (\ref{equa}) is less than $1$. Hence there is a straightforward upper bound on $\P_{\theta_j} (  t_j > T)$:
\begin{eqnarray}
  \P_{\theta_j} (   t_j > T)& \leq & \exp(-(hT -\phi_{\theta_j}(h) )). \label{upper-GD}
\end{eqnarray}
To complete this part of the proof, it remains to determine the minimum value of a positive value $h$ on the right-hand side of (\ref{upper-GD}).
The minimum is  attained  for  positive $h$ such that
  \begin{eqnarray} \E_{\theta_j,h} ( t_j)&=& T \label{sol-h-Prop61-UPP} \end{eqnarray}
since
$$\left\{ \begin{array}{lcl}
 (\phi_{\theta_j} (h)-hT)'&=& 
\E_{\theta_j,h} ( t_j) -T,\\
   (\phi_{\theta_j} (h)-hT)^{''} & =&\V_{\theta_j,h} (t_j)\; \geq 0 ,
\end{array} \right. $$
where $( \cdot)^{'}$ and $( \cdot)^{''}$ denote the first and second derivatives with respect to $h$, respectively, and, $\E_{\theta_j,h}$ and $\V_{\theta_j,h}$ are the expectation and variance with respect to
$\P_{\theta_j,h}$.

In order to find  $h$ that solves equation (\ref{sol-h-Prop61-UPP}), we need to determine  $\phi_{\theta_j}$.
For this, set $ \nu_{j,k}=\displaystyle{\frac{\theta_{j,k}}{\epsilon}}$. Then for any positive  $h$ such that
$h  \rightarrow  +\infty$ and  $h\displaystyle{\max_k}w_k =o(1)$, we obtain
\begin{eqnarray}
\phi_{\theta_j} (h) &=& \log \prod_{k} \E_{\theta_j} [\exp (h  w_k ((\nu_{j,k}+ \eta_{j,k})^2-1))] \nonumber \\
& =& \sum_k \{-h w_k + \frac{h w_k \nu_{j,k}^2}{(1-2h w_k)} -\frac 1 2 \log (1-2h w_k) \} \nonumber \\
& =&  \sum_k \{-h w_k + h w_k \nu_{j,k}^2 (1 +2 h w_k + o(h w_k))   \nonumber \\
&& \;\;\;\;\;\;\;\;\;-\frac 1 2 (-2 hw_k -\frac{(2h w_k)^2}{2} +o(h^2w_k^2))\} \nonumber \\
& =& \sum_k \{ hw_k \nu_{j,k}^2 (1 +o(h\max_k w_k))+ h^2 w_k^2 (2  \nu_{j,k}^2+ 1) +  o(h^2w_k^2) \}\nonumber \\
& =& h \E_{\theta_j} (t_j) (1 +o(h\max_k w_k))  + \frac{h^2}{2}(1 +  o(1)) + o(h^2) \label{gdf},
\end{eqnarray}
where the last equality sign in (\ref{gdf})  follows  from $(T -\E_{\theta_j} (t_j))\rightarrow +\infty$  and  $T \; {\displaystyle \max_k} w_k  =o(1)$  as $T \rightarrow +\infty$.
Next,  differentiating  the right-hand side of (\ref{gdf}) with respect to $h$ yields
$$\begin{array}{rclccc}
(\phi_{\theta_j} (h) -hT)'&=& 0
&\Rightarrow \;& h \sim& T - \E_{\theta_j} (t_j),\; \; \mbox{ {\rm as $T - \E_{\theta_j} (t_j)$ goes to infinity.} }
\end{array}$$
As $(T- \E_{\theta_j} (t_j))\stackrel{T \rightarrow + \infty}{ \; \; \; \; \; \longrightarrow + \infty}$, this leads  to the following optimal upper bound for
right-hand side of (\ref{upper-GD}):
\begin{eqnarray*}
\exp\left((T \! - \! \E_{\theta_j} (t_j)) \E_{\theta_j} (t_j) \! +  \!\frac{(T \! - \! \E_{\theta_j} (t_j))^2}{2}  \!- \! T (T \!- \!\E_{\theta_j} (t_j)\right)  \sim  \exp(
- \frac{(T\! - \!\E_{\theta_j} (t_j))^2}{2}). \label{opti_upp_GD}
\end{eqnarray*}
Since by assumption $T\; \displaystyle{\max_k}w_k =o(1)$, the condition $h\; \displaystyle{\max_k}w_k =o(1)$ with $(T-\E_{\theta_j}(t_j))$ in place of $h$ is  fulfilled. \\
By assumption $T\; \displaystyle{\max_k}w_k =o(1)$, hence  the optimal upper bound under $\P_0$ is $\exp (-\frac{T^2}{2})$ as $T$ goes to infinity.
This completes the proof of the upper bound. \\


\noindent
\vspace{0.3cm}
{\it Lower Bound.} We are interested in obtaining a lower  bound for  (\ref{equa}).
This is done by  first considering a new family of probability distributions under which the  normalized statistics $t_j$  are  proved to be asymptotically Gaussian. \\

For $h>0$ satisfying equation (\ref{sol-h-Prop61-UPP}),
let us introduce the following probability  measures $\P_{\theta_j,h,k}$: $$\displaystyle{\frac{d \P_{\theta_j,h,k}}{d \P_0}}=\exp( h  t_{j,k} ) \;\exp(-\phi_{\theta_{j,k}}(h)),$$
with $t_{j,k}$  defined in (\ref{tjk}),  $\phi_{\theta_{j,k}}(h) = \log \E_{\theta_{j,k}} (\exp( h  t_{j,k}))$ and where
$\E_{\theta_{j,k}}$  stands for the expectation with respect to the observations  $(x_{j,k})_{j,k}$ of
(\ref{transformed_obs}). Denote by $\E_{\theta_j,h,k}$ and $\V_{\theta_j,h,k}$  the expectation and variance with respect to
$\P_{\theta_j,h,k}$.

To establish the asymptotic normality of $t_j$, we will check that the Lyapunov condition is satisfied.
To this end, set $\sigma^2_{j,h,k} = \V_{\theta_j,h,k}(t_{j,k})$  and  $\sigma^2_{j,h} =\sum_k \sigma^2_{j,h,k}$. \\
Denote by $\phi^{(2)}_{\theta_{j,k}}$ and $\phi^{(4)}_{\theta_{j,k}}$ the second and  fourth
derivatives  of $\phi_{\theta_{j,k}}$  with respect to $h$, respectively. Using  well-known relations between moments of $t_j$ under  $\P_{\theta_j,h,k}$ and the successive derivatives of $\phi_{\theta_{j,k}}(h)$ with respect to $h$, in particular, $\sigma^2_{j,h}=\sum_k \phi^{(2)}_{\theta_{j,k}}$,
we get
\begin{eqnarray*}
   \frac{\sum_k \E_{\theta_j,h,k}( t_{j,k} - \E_{\theta_j,h,k}(t_{j,k}))^4}{ (\sum_k\sigma^2_{j,h,k})^2}
& =& \frac{ 3 \sum_k (\phi^{(2)}_{\theta_{j,k}}(h))^2 + \sum_k \phi^{(4)}_{\theta_{j,k}}(h)}{(\sum_k\phi^{(2)}_{\theta_{j,k}}(h))^2} \nonumber  \\
& \leq & \frac{ 4 \max (w_k^2) \sum_k w_k^2 (1+o(1))+ o(1)}{1}  \nonumber  \\
& =& o(1),
\end{eqnarray*}  where the last relation follows from $\max w_k=o(1)$
 and  relation (\ref{gdf}), since by  (\ref{gdf}) we get $\sum_k \phi^{(4)}_{\theta_{j,k}}(h) =
\phi^{(4)}_{\theta_j} (h)=o(1)$.
The  Lyapunov condition is then satisfied.   This implies that under
$\P_{\theta_j,h}$, $Z_{j,h}=\displaystyle{\frac{t_j - \E_{\theta_j,h} (t_j)}{\sigma_{j,h}}}$ is asymptotically a real standard Gaussian random variable.\\

Let us return to relation (\ref{equa}), where
$h$ is chosen to have  $\E_{\theta_j,h} (t_j)=T$, and observe that $$\E_{\theta_j,h} [ \I_{(t_j > T)} \exp(-h( t_j- T))] = \E_{\theta_j,h} [ \I_{(Z_{j,h} >0)} \exp(-h Z_{j,h}\sigma_{j,h})].$$
Due to the asymptotic normality of $t_j$, for any $\delta>0$,
\begin{eqnarray}
 \!\!\!\!\!\!\! \E_{\theta_j,h} [ \I_{(Z_{j,h} >0)} \exp(-h Z_{j,h}\sigma_{j,h})]
&\!\!\!\!\!  =& \!\!\!\!\E_{\theta_j,h} [ \I_{(Z_{j,h}  \in (0, \delta))} \exp(-h Z_{j,h}\sigma_{j,h})] + \nonumber \\
&& \E_{\theta_j,h} [ \I_{(Z_{j,h}  >  \delta)} \exp(-h Z_{j,h}\sigma_{j,h})] \nonumber \\
&\!\!\!\!\! >& \!\!\!\!(\P_{\theta_j,h} ( Z_{j,h}\in (0, \delta)) +o(1))   \exp(-h \delta \sigma_{j,h}) . \label{equa2}
\end{eqnarray}
By choosing $\delta = o(h)$ in relation (\ref{equa2})  implies that
\begin{eqnarray}
\log (\P_{\theta_j} (t_j > T) ) & \geq & \phi_{\theta_j} (h) - hT - o(h^2). \label{lower-GD}
\end{eqnarray}
Up to $o(h^{2})$, the right-hand side of (\ref{lower-GD}) corresponds  to the argument of the exponential function on the right-hand side of (\ref{upper-GD}). This entails that the right-hand side of (\ref{lower-GD}) is   equivalent to $\displaystyle{ -\frac{(T - \E_{\theta_j} (t_j) )^{2}}{2}}$.
This completes the  proof of the lower bound, and thus  Proposition \ref{GD} is proved.

\begin{lemma} \label{resu_tech}   \begin{itemize}
                                   \item (i) Expectation and variance of $t_j$ defined by (\ref{tj}).
\begin{eqnarray}
\E_{{ \theta_j}} (t_j) &=&  \xi_j \epsilon^{-2} \kappa({\theta_j},w),  \label{Lemm_E_1} \\
\V_{{ \theta_j}} (t_j) & =&  1 + O ((\max_{k \in \BBz} \; w_k) \; \E_{{ \theta_j}} (t_j)) .  \label{Lemm_V_1_1}
\end{eqnarray}
\item (ii) Expectation and variance of $L(u)$ defined by (\ref{Inter-Hig}).
Assume that $T_d \max w_k=o(1)$
and consider any ${\overline{\boldsymbol{\theta}}}=(\xi_1 \theta_1, \ldots,\xi_d\theta_d)$ such that $\sum_{j=1}^{d}\xi_j=d^{1-b}$.
Moreover, assume that  for all nonzero $\xi_j$,     $\E_{\theta_j} (t_j) \ge c T_d $, with some positive $c$, and
 $\displaystyle{\max_{j:\xi_j=1}}\E_{\theta_j} (t_j)=O(T_d)$. Then, for all $u \in (0, \sqrt 2]$,
\begin{eqnarray*}
\E_{\overline{ \boldsymbol{\theta}}}(L(u)) & \ge & d^{\frac{1}{2} -b +
(\frac{u^2}{4}-\frac{((u- c)_+)^2}{2} )(1+o(1))}(1+o(1)) ,
  \label{E_1IH} \\
  \V_{\overline{\boldsymbol{\theta}}}(L(u))&  =& o(d^{\eta}\,\E_{\overline{ \boldsymbol{\theta}}}(L(u)))  + o(1),
\quad \eta=o(1),
  \label{V_E_1IH}
  \end{eqnarray*}
where $x_+=\max(0,x)$.
\end{itemize}
\end{lemma}
\begin{remark}\label{RN}
Under $\P_0$, the statistics $t_j$ and  $L(u)$ have  zero mean and
unit variance. Moreover, under $\P_0$ and the assumption $\displaystyle{\max_k w_k}=o(1)$, the
statistics $t_j$ are asymptotically standard Gaussian. Under
$\P_{\theta_j}$, the statistics $t_j-\E_{\theta_j}t_j$ are
asymptotically standard Gaussian if $\displaystyle{\max_k w_k} \,
\E_{\theta_j}t_j=o(1)$, see Lemma 3.1 in \cite{IS.02a}.
\end{remark}


\vspace{0.2cm}
\noindent{\it Proof of Lemma \ref{resu_tech}.}\\
(i) Recall that
$\displaystyle{\sum_k} w_k^2=1/ 2$.
For each index $j$ satisfying $\xi_j=1$,  the random variable $\displaystyle{(\frac{x_{j,k}}{\epsilon})^2}$ is
a $\P_{{ \theta_j}}$-noncentral $\chi^2(1, \theta^2_{j,k} \epsilon^{-2})$.
From this relation (\ref{Lemm_E_1}) is   easily  obtained.
 Relation (\ref{Lemm_V_1_1})  is deduced from the following calculations:
\begin{eqnarray}
\V_{{ \theta_j}} (t_j)& =& \sum_{k \in \BBz} w_k^2( 2 + 4 \epsilon^{-2} \xi_j \theta_{j,k}^2 ) \nonumber \\
& =&  1 + \sum_{k \in \BBz} w_k^2   4 \epsilon^{-2} \xi_j  \theta_{j,k}^2 \nonumber\\
& = & 1 + O (\displaystyle{\max_{k \in \BBz} w_k} \;  \epsilon^{-2}  \xi_j \kappa({ \theta_j},w))   \nonumber \\
& = & 1 + O(\displaystyle{\max_{k \in \BBz} w_k} \; \E_{{\theta_j}} (t_j)) . \nonumber
\end{eqnarray}
(ii) For any $u\in(0,\sqrt 2]$, as $T_d
\rightarrow +\infty$, Proposition \ref{GD} gives a control over $C_u$ defined by (\ref{Cu}):
\begin{eqnarray*}
C_u^2 &=
& d^{-1} \exp(\frac{u^2 T_d^2}{2}(1+o(1)))(1 -\exp(\frac{-u^2 T_d^2}{2}(1+o(1))))^{-1} \nonumber \\
& = & d^{-1 + \frac{u^2}{2}(1+o(1))}.
\end{eqnarray*}
Since $u \leq \sqrt 2$, the exponent of $d$ in $C_u$ is $o(1)$.
%
%
 \\
{\bf Case 1:}  for the nonzero $\xi_j$'s, assume that
$\lim\sup(u T_d-\E_{\theta_j} (t_j))< + \infty$. 
In this case, the probability $\P_{\theta_j} (t_j>u T_d)=
\P_{\theta_j} (t_j- \E_{\theta_j} (t_j) > u T_d -\E_{\theta_j} (t_j)
)$ is bounded away from zero. This follows from the asymptotic
 normality of $t_j- \E_{\theta_j}(t_j)$ for $\E_{\theta_j}(t_j)=O(T_d)$ (see Remark
 \ref{RN})


\noindent {\bf Case 2:} for the nonzero $\xi_j$'s, assume
that  $uT_d-\E_{\theta_j} (t_j)  \to + \infty$. Then, for any nonzero
$\xi_j$,
  Proposition \ref{GD} implies that  
\begin{eqnarray*}
 \log \P_{\theta_j} (t_j >uT_d) & \ge &  - \frac{(uT_d-c T_d)^2}{2}(1 + o (1))  . \label{GD-tot}
\end{eqnarray*}

\bigskip
Recall  that the number of nonzero $\xi_j$ is equal to $K=d^{1-b}$ and that  for all nonzero $\xi_j$,    $\E_{\theta_j} (t_j) \ge c T_d $ for some positive $c$ such that
 $\max_{j:\xi_j=1}\E_{\theta_j} (t_j)=O(T_d)$.
To sum up,  the cases 1 and 2  
entail that
\begin{eqnarray}
\E_{\overline{\boldsymbol{\theta}}}(L(u)) & =& C_u \sum_{j: \xi_j =1} \left(\P_{\theta_j} (t_j > u T_d ) - {\tilde \Phi}_0(u T_d)\right) \nonumber\\
& \geq & C_u K \left(d^{-\frac{((u -c)_+)^2}{2}(1+o(1))} -d^{-\frac{u^2 }{2}(1+o(1))}\right) \nonumber \\
& = & d^{-\frac{1}{2}  + \frac{u^2}{4}(1+o(1)) +1-b} \left( d^{-\frac{((u-c)_+)^2}{2} (1+o(1))} -d^{-\frac{u^2}{2}(1+o(1))}     \right) (1+o(1))  \nonumber \\
& =&  d^{\frac{1}{2} -b + (\frac{u^2}{4}-\frac{((u-c)_+)^2}{2} )(1+o(1))}(1+o(1)).\nonumber
\end{eqnarray}
Similarly, let us study the variance of $L(u)$.  Using Proposition
\ref{GD}, we obtain
\begin{eqnarray*}
\V_{\overline{\boldsymbol{\theta}}}(L(u))&  =&  C_u^2 \sum_{j: \xi_j =1} \P_{\theta_j} (t_j > u T_d )
\P_{\theta_j} (t_j \leq u T_d ) + C_u^2 \sum_{j: \xi_j =0}  {\tilde \Phi}_0(u T_d) (1-{\tilde \Phi}_0(u T_d))
\nonumber \\
  & =& C_u^2 K \P_{\theta_j} (t_j > u T_d )(1+o(1))  +(d^{b-1}  + d^{-b})(1+o(1))
  \nonumber  \\
  & =& ( C_u \E_{\overline{\boldsymbol{\theta}}}(L(u)) + d^{b-1})(1 + o(1))  \nonumber  \\
  & =& o(d^\eta \E_{\overline{\boldsymbol{\theta}}}(L(u))) +o(1),\quad \eta=o(1).
\end{eqnarray*}

\subsection{Upper bound}
\begin{remark}\label{RNN}
Note that the condition $T_d\;  {\displaystyle \max_k} w_k(r_\e^{\star})= o(1)$
follows from assumption   $\log d=o(\epsilon^{-2/(2 \tau +1)})$. Indeed,   Remark \ref{rmv0}
and relations (\ref{Lagrange}) imply that $T_d \max w_k (r_\e^{\star})\leq
(r_\epsilon^{\star})^{1/(2 \tau)} T_d$,
where the term on the right-hand side goes to zero as soon as
$\log d=o(\epsilon^{-2/(2 \tau +1)})$. Therefore, assumption $\log
d=o(\epsilon^{-2/(2 \tau +1)})$ allows us to apply Proposition
\ref{GD} and Lemma  \ref{resu_tech}.
\end{remark}

\noindent
{\it Proof of  {\it (ii)}--Theorem \ref{Moderate_case}.} \\
{\bf Type I error probability of $\psi^{\chi^2}_\alpha$.}
It follows from the Central Limit Theorem that, under the null hypothesis, $t_b$ is asymptotically a standard normal random variable. Therefore
\begin{eqnarray*}
\P_0( t_b > T_\alpha) & =& \Phi(-T_\alpha)+ o(1) =\alpha + o(1).
\end{eqnarray*}

 \vspace{0.1cm} \noindent
{\bf Type II error probability of $\psi^{\chi^2}_\alpha$ uniformly
over  $\Theta_d(\tau,r_\epsilon,b)$} for $r_\e\ge Br^\star_\e,\ B\ge 1$.
Thanks to  Lemma \ref{resu_tech}, uniformly over
${\overline{\boldsymbol{\theta}}}\in \Theta_d(\tau,r_\epsilon,b)$,
we have
\begin{eqnarray*}
\Var_{\overline{\boldsymbol{\theta}}} (t_b)& = &
\frac 1 d  \sum_{j=1}^d   (1 + O (\E_{\theta_j} (t_{j,b}))) \\
\E_{\overline{\boldsymbol{\theta}}} (t_b) & =&  d^{-1/2} \sum_{ j=1}^d  \E_{\theta_j} (t_{j,b}).
\end{eqnarray*}
This implies that $\Var_{\overline{\boldsymbol{\theta}}} ( t_b)= o
((\E_{\overline{\boldsymbol{\theta}}} (t_b))^2)$ provided that
$\E_{\overline{\boldsymbol{\theta}}} (t_b)  \rightarrow +\infty$. Let
us study $\E_{\overline{\boldsymbol{\theta}}} (t_b)$: from Proposition
\ref{P1}, Lemma \ref{resu_tech}, and  relation (\ref{PB_minim}), we
get uniformly over $\Theta_d(\tau,r_\epsilon,b)$ with $r_\e\ge
Br^\star_\e,\ B\ge 1$:
\begin{eqnarray}
\E_{\overline{\boldsymbol{\theta}}} (t_b) \geq  d^{1/2 -b}  B^2
a(r_\epsilon^\star) \rightarrow +\infty \mbox{ {\rm as soon as } }
B^2d^{1/2-b} a(r_\epsilon^\star)\asymp B^2 \rightarrow +\infty,  \mbox{ {\rm i.e., as soon as } } r_\e/r^\star_\e \rightarrow
+\infty,\quad \label{TT}
\end{eqnarray}
where $r_\epsilon^{\star}\asymp (\epsilon^4 d^{2b -1})^{\tau/(4 \tau +1)}$.

Due to (\ref{TT}),  using Markov's inequality and Lemma \ref{resu_tech}, for all  $\overline{\boldsymbol{\theta}}$ in $\Theta_d(\tau,r_\epsilon,b)$,
\begin{eqnarray*}
\P_{\overline{\boldsymbol{\theta}}} ( t_b \leq T_\alpha) & =&
\P_{\overline{\boldsymbol{\theta}}}  ( t_b -\E_{\overline{\boldsymbol{\theta}}} (t_b) \leq  T_\alpha -\E_{\overline{\boldsymbol{\theta}}} (t_b))
\nonumber \\
& \leq &  \P_{\overline{\boldsymbol{\theta}}}  ( |t_b -\E_{\overline{\boldsymbol{\theta}}} (t_b)| \geq  \E_{\overline{\boldsymbol{\theta}}} (t_b)- T_\alpha  )  \nonumber \\
& \leq & \frac{ \Var_{\overline{\boldsymbol{\theta}}}(
t_b)}{(\E_{\overline{\boldsymbol{\theta}}} (t_b)- T_\alpha )^2}=o(1).
\end{eqnarray*}
This entails that
 $\beta(\psi^{\chi^2}_\alpha,\Theta_d(\tau,r_\epsilon,b))$ goes to zero as soon as  $d^{1/2-b} a(r_\epsilon) \rightarrow +\infty$, i.e.,
as soon as $\displaystyle{\frac{a(r_\epsilon)}{a( r_\epsilon^{\star})}} \rightarrow +\infty$  where
$a(r_\epsilon^{\star}) \asymp  d^{b -1/2}$.

\vspace{0.3cm}
\noindent
{\it Proof of  {\it (ii)}--Theorem \ref{Intermediate_case}.} \\
{\bf Type I error probability of $\psi^{HC}$.} Observe that $
w(\psi^{HC})\le w(\psi^{L})+w(\psi^{max}).$ The assumption 
$\log(d) = o (\epsilon^{-2/(2\tau +1)})$ 
implies that $T_d \max_k w_k(r_\epsilon^{\star})=o(1)$. Therefore the application of  Proposition
\ref{GD} and the fact that $D^2>2$ and $N=o(T_d)$ yield
\begin{eqnarray*}
w(\psi^{max})&=&\P_0(\max_{1\le j\le d} \max_{1\le l\le
N}t_{j,b_l}>DT_d)\le
\sum_{j=1}^d\sum_{l=1}^N\P_0(t_{j,b_l}>DT_d)\\
&\le &Nd\exp(-D^2T_d^2/2(1+o_d(1)))= Nd^{1-D^2/2(1+o_d(1))}\to 0.
\end{eqnarray*}
By Lemma \ref{resu_tech} and applying Markov's inequality,
\begin{eqnarray*}
w(\psi^{L})= \P_0(\max_{1 \leq l \leq N-1}L(u_l,b_l) > H) & \leq &
\sum_{l=1}^{N-1} \P_0(L(u_l,b_l) >H)
 \nonumber \\
& \leq &  \sum_{l=1}^{N-1} \frac{\V_0 (L(u_l,b_l)) }{H^2} \nonumber \\
& \leq & \frac{(N-1)}{H^2}, \label{Inter_Hig_typeI}
\end{eqnarray*}
which goes to zero as $d\rightarrow + \infty$ since $H \sim (\log d)^C$, with $C > \frac 1 4$  and $N =o_d(T_d)$.

\vspace{0.3cm} \noindent {\bf Type II error probability of
$\psi^{HC}$ uniformly over $\Theta_d(\tau,r_\epsilon,b)$.} For any
$\overline{\boldsymbol{\theta}} \in \Theta_d(\tau,r_\epsilon,b)$, we
obtain
 \begin{eqnarray}\label{test1}
\E_{\overline{\boldsymbol{\theta}}}(1-\psi^{HC})&\le&
\min(\E_{\overline{\boldsymbol{\theta}}}(1-\psi^{max}),\E_{\overline{\boldsymbol{\theta}}}(1-\psi^{L})),\\
\E_{\overline{\boldsymbol{\theta}}}(1-\psi^{max})&\le& \min_{j:
\xi_j=1}\min_{1\le l\le
N}\P_{{\theta_j}}(t_{j,b_l}\le DT_d).
\label{test2}
\end{eqnarray}

First, let us consider the alternatives $\overline{\boldsymbol{\theta}} \in
\Theta_d(\tau,r_\epsilon,b)$ such that for a nonzero $\xi_j$, there exists $l \in \{1,\ldots,N\}$
 for which  $\E_{{\theta_j}}t_{j,b_l}\ge
D_1T_d$ with $D_1>D$. From Lemma \ref{resu_tech}(i) and
Markov's inequality, we obtain
 \begin{eqnarray}
\P_{{\theta_j}}(t_{j,b_l}\le DT_d)&\le&
\P_{{\theta_j}}(| t_{j,b_l}-\E_{{\theta_j}}(t_{j,b_l})|\ge
\E_{{\theta_j}}(t_{j,b_l})-DT_d) \nonumber \\
&\le& \frac{ \Var_{{\theta_j}}({
t}_{j,b_l})}{(\E_{{\theta_j}} (t_{j,b_l})- DT_d
)^2}=o(1). \label{test3}
\end{eqnarray}

Second,  in view of \nref{test1}, \nref{test2}, \nref{test3},  it suffices
to study the test procedures $\psi^{L}$ under the alternatives $\overline{\boldsymbol{\theta}}
\in \Theta_d(\tau,r_\epsilon,b)$ such that $\displaystyle{\max_{j:
\xi_j=1}\max_{1\le l\le
N}\E_{\theta_j}t_{j,b_l}=O(T_d)}$. Then we obtain
\begin{eqnarray*}
\E_{\overline{\boldsymbol{\theta}}}(1-\psi^{L})&=&\P_{\overline{\boldsymbol{\theta}}}(\max_{1
\leq l \leq N-1} L(u_l,b_l) \leq  H) 
\leq  \min_{{1 \leq l
\leq N-1}} \P_{\overline{\boldsymbol{\theta}}}(L(u_l,b_l)  \leq H).
\end{eqnarray*}

For any $l \in
\{1,\ldots,N-1\}$,
 \begin{eqnarray}
  \P_{\overline{\boldsymbol{\theta}}}(L(u_l,b_l)  \leq H)&\le &\P_{\overline{\boldsymbol{\theta}}}(L(u_l,b_l)
-\E_{\overline{\boldsymbol{\theta}}} (L(u_l,b_l)) \leq H-
\E_{\overline{\boldsymbol{\theta}}} (L(u_l,b_l))) \nonumber \\
& \leq & \P_{\overline{\boldsymbol{\theta}}}(-|L(u_l,b_l) -
\E_{\overline{\boldsymbol{\theta}}} (L(u_l,b_l))| \leq H-\E_{\overline{\boldsymbol{\theta}}} (L(u_l,b_l))) \nonumber \\
& \leq & \P_{\overline{\boldsymbol{\theta}}}(|L(u_l,b_l) -\E_{\overline{\boldsymbol{\theta}}} (L(u_l,b_l)) |\geq -H+
\E_{\overline{\boldsymbol{\theta}}} (L(u_l,b_l))) \nonumber \\
& \leq &
\frac{\V_{\overline{\boldsymbol{\theta}}}(L(u_l,b_l))}{(\E_{\overline{\boldsymbol{\theta}}}
(L(u_l,b_l)) -H)^2}. \label{Inter_Hig_typeII}
\end{eqnarray}
 For any  $b_l \in (1/2, 1)$, if we prove that
$\displaystyle{\inf_{\overline{\boldsymbol{\theta}} \in \Theta_d(\tau,r_\epsilon,b)}}
\E_{\overline{\boldsymbol{\theta}}} (L(u_l,b_l))$ goes  to infinity
as a power of $d$  ($ d \rightarrow +\infty$), then Lemma
\ref{resu_tech} and the choice of $H$ (recall $H=O_d((\log d)^C)$,
with $C
>1/4$) yield the result since in this case the right-hand side of relation
(\ref{Inter_Hig_typeII}) goes to zero.

Third, for $b\in (1/2,1)$,  take an index $l$ in $\{1,\ldots,N-1\}$ such that $b_l\le b\le b_{l+1}$.
This, combined with the continuity of $\varphi$, yields
$$
b_l=b+o(1),\quad r_\e^\star(b_l)\le r_\e^\star(b)\sim r_\e^\star(b_l),\quad
a(r_\e^\star(b_l))\le a(r_\e^\star(b))\sim a(r_\e^\star(b_l)).
$$



Let $\overline{\boldsymbol{\theta}} \in \Theta_d(\tau,r_\epsilon,b)$ with $b\in (1/2,1)$ and
$\lim\inf(a(r_\e)/a(r_\e^\star(b))>1$. Then $r_\e\ge (1+\delta)r_\e^\star(b_l)$
for some $\delta>0$.  Proposition \ref{P1}  entails that for $j$
such that $\xi_j=1$ we have
$$
\E_{\theta_j}t_{j,b_l}\ge (1+\delta)^2a(r^\star_\e(b_l))\sim
(1+\delta)^2a(r^\star_\e(b))\sim (1+\delta)^2\varphi(b)T_d.
$$

We then derive from  Lemma \ref{resu_tech} with
$c=c(b)=(1+\delta)^2\varphi(b)$ that
\begin{eqnarray}
\inf_{\overline{\boldsymbol{\theta}} \in
\Theta_d(\tau,r_\epsilon,b)} \E_{\overline{\boldsymbol{\theta}}}
(L(u_l,b_l)) &> &d^{\frac 1 2 + \frac{u_l^2}{4}  -b - \frac{((u_l -
c(b))_+)^2}{2}(1+o(1)) } (1+o(1))\label{need-ext}.
\end{eqnarray}
Finally, denote the main term in the exponent of $d$ in
(\ref{need-ext}) by $$
M=\frac 1 2 + \frac{u(b)^2}{4}  -b - \frac{((u(b) -
c(b))_+)^2}{2}.
$$
 To obtain the result, it is sufficient  to
prove that  $M$ is positive and bounded away from zero for any
$\delta>0$.

\noindent {\it Intermediate sparsity case.} This case
corresponds to $b \in (1/2, 3/4]$. Recall that  $u(b)= 2 \varphi_1(b)$,
where $\varphi_1$ is defined in (\ref{phi}). Then

 $M>0 \Leftrightarrow \left\{ \begin{array}{lll}
  (\varphi_1^2(b)/2)((1+\delta)^2-1)(3-(1+\delta)^{2})>0  & \mbox{{\rm for }}  & 0<\delta< \sqrt 2-1 \\
  \varphi_1^2(b)/2>0& \mbox{{\rm for }}
 & \delta \ge \sqrt 2-1\end{array} \right.. $

\noindent  The latter inequalities are obviously satisfied. This leads to the result.

\noindent {\it Highest sparsity  case.} In this case $b \in (3/4, 1)$ and $u(b)= \sqrt 2$. Then

$M>0 \Leftrightarrow \left\{ \begin{array}{lll}
  ((1+\delta)^2-1)\varphi_2(b)>0  & \mbox{{\rm for }}  & (1+\delta)^2\varphi_2(b)\le \sqrt 2 \\
1-b>0  & \mbox{{\rm for }}
 &(1+\delta)^2\varphi_2(b)>
\sqrt 2\end{array}  \right.. $

\noindent Again,  the latter inequalities are  satisfied, and the result follows.

\vspace{0.3cm}
\noindent
{\it Proof of  {\it (i)}--Theorem \ref{thm:extended-problem}.}  \\
Similar to the
proof of part (ii) of Theorem \ref{Moderate_case},    due to
(\ref{kappa-ext}) and \nref{kappa-ext_1}, uniformly over
$\overline{\boldsymbol{\theta}}\in \Theta^{ext}_d(\tau,K^{1/2},
K^{1/2}r_\epsilon,b)$, the type II error probability of
$\psi^{\chi^2}_{\alpha}$ goes to zero as soon as
$r_\e/r^\star_\e\to +\infty$.

\vspace{0.3cm}
\noindent
{\it Proof of  {\it (ii)}--Theorem \ref{thm:extended-problem}.}  \\
The proof of the fact that the type II error probability
of $\psi^{HC}$ goes to zero as $d \rightarrow +\infty$ is similar to
the one of Theorem \ref{Intermediate_case}. Recall that $K=d^{1-b}$
is the number of nonzero $\xi_j$'s and suppose without loss of
generality that $\xi_j=1, \forall j \in \{1,\ldots,K\}$ and $\xi_j
=0, \forall j \in \{K+1,\ldots,d\}$. Note that relations (\ref{test1}) and (\ref{test2}) remain valid for any
$\overline{\boldsymbol{\theta}} \in
\Theta^{ext}_d(\tau,K^{1/2}, K^{1/2}r_\epsilon,b)$.

First, similarly to (\ref{test3}),  for any  $\overline{\boldsymbol{\theta}} \in
\Theta^{ext}_d(\tau,K^{1/2}, K^{1/2}r_\epsilon,b)$ such that for the  nonzero $\xi_j$'s, there exists $l \in \{1,\ldots,N\}$
 for which  $\E_{{\theta_j}}t_{j,b_l}\ge
D_1T_d$ with $D_1>D$, the type II error probability of  $\psi^{HC}$ vanishes asymptotically.
Therefore,   it suffices
to study the test procedures $\psi^{L}$ under the alternatives $\overline{\boldsymbol{\theta}}
\in\Theta^{ext}_d(\tau,K^{1/2}, K^{1/2}r_\epsilon,b)$ such that $\displaystyle{\max_{j:
\xi_j=1}\max_{1\le l\le
N}\E_{\theta_j}t_{j,b_l}=O(T_d)}$. Therefore, let us take $\delta >0$
 and consider the alternatives that are as far away from the null hypothesis as $r_\e$ such that
$r_\e\ge (1+\delta)r^\star_\e(b)$, where $r^\star_\e(b)$ is determined by
$a(r_\e^\star(b))\sim T_d\varphi(b)$.

Second, for any $l \in \{1,\ldots,N\}$,  observe that the only difference between the proofs of  the extended and
initial problems lies in the study of
\begin{eqnarray}
\inf_{\overline{\boldsymbol{\theta}} \in
\Theta^{ext}_d(\tau,K^{1/2}, K^{1/2}r_\epsilon,b)} \sum_{j=1}^K
\P_{\theta_j}  (t_{j,b_l} - \E_{\theta_j}(t_{j,b_l})> u_l T_d -
\E_{\theta_j}(t_{j,b_l})).\label{Terme-traiter}
\end{eqnarray}
Now it is no more  possible to control (\ref{Terme-traiter}) by using  Lemma \ref{resu_tech} (ii)
because the condition  $\E_{\theta_j} (t_j) \ge c T_d$ is not necessarily satisfied for all nonzero $\xi_j$'s.
In fact, the only condition we have is   $\displaystyle{\sum_{j=1}^{K}\E_{\theta_j} (t_j)} \ge c K T_d$ with some  constant  $c>1$.

Let us now explain why the current proof is reduced to the study of (\ref{Terme-traiter}).
As in  (\ref{Inter_Hig_typeII}), we get for any $\overline{\boldsymbol{\theta}}$  in
$ \Theta^{ext}_d(\tau,K^{1/2}, K^{1/2}r_\epsilon,b)$,
  \begin{eqnarray*}
\P_{\overline{\boldsymbol{\theta}}}(\max_{1 \leq l \leq N}
L(u_l,b_l) \leq  H) & \leq & \min_{1\le l\le N}
\frac{\V_{\overline{\boldsymbol{\theta}}}
(L(u_l,b_l))}{(\E_{\overline{\boldsymbol{\theta}}} (L(u_l,b_l))
-H)^2}. \label{Ext-Inter_Hig_typeII}
\end{eqnarray*}
Due to Lemma \ref{resu_tech} and the fact that $H=O_d((\log d)^C)$ with $C >
1/4$, in order to obtain the result, it remains to prove that for any $l$ such that $b_l\le b\le
b_{l+1}$,  $\displaystyle{\inf_{\overline{\boldsymbol{\theta}}
\in \Theta^{ext}_d(\tau,K^{1/2}, K^{1/2}r_\epsilon,b)}}
\E_{\overline{\boldsymbol{\theta}}} (L(u_l,b_l)) \stackrel{d
\rightarrow +\infty}{\longrightarrow} +\infty$ as a positive power
of $d$. Finally, recall that
\begin{eqnarray}
\E_{\overline{\boldsymbol{\theta}}} (L(u_l,b_l))
& =&  C_{u_l,b_l} \sum_{j=1}^K \left( \P_{\theta_j} (t_{j,b_l} -
\E_{\theta_j}(t_{j,b_l}) > u_l T_d- \E_{\theta_j}(t_{j,b_l})) -
{\tilde \Phi}_{0,b_l} (u_l T_d) \right),\label{equa*}
\end{eqnarray}
where $C_{u_l,b_l}= (d{\tilde \Phi}_{0,b_l} (u_l T_d) (1-{\tilde \Phi}_{0,b_l}
(u_l T_d)))^{1/2}$ and $ {\tilde \Phi}_{0,b_l} (x) = \P_0 (t_{j,b_l} >
x)$.  The term on the right-hand side of (\ref{equa*})
corresponds to the product of (\ref{Terme-traiter}) and   $C_{u_l,b_l}$.
 The quantity $C_{u_l,b_l}$ is controlled   by Lemma
\ref{resu_tech} and Proposition \ref{GD}. Thus it remains to study (\ref{Terme-traiter}).

{Third, the application of Proposition \ref{GD} gives the following approximation of (\ref{Terme-traiter}),
\begin{eqnarray*}
\sum_{j=1}^K \P_{\theta_j}  (t_{j,b_l} - \E_{\theta_j}(t_{j,b_l})>
u_l T_d - \E_{\theta_j}(t_{j,b_l})) &=& \sum_{j=1}^K \exp(-\frac{((
u_l T_d - \E_{\theta_j}(t_{j,b_l}))_+)^2}{2})O(1).
\label{Terme-traiter-1}
\end{eqnarray*}
Recall that $a(r_\epsilon)$  given by (\ref{sol_a}) is the solution
of the extremal problem (\ref{minimisation}).   Set
$\eta_j=\E_{\theta_j}(t_{j,b_l})$,
$\eta_0=(1+\delta)^2a(r_\epsilon^\star(b))\sim
(1+\delta)^2a(r_\epsilon^\star(b_l))$ and
$f_T(\eta)=\exp(-\frac{(T-\eta)^2}{2})$ $ \forall \eta \in [0, R]$,
where  $R=R(T)>0$  will be specified later on. Consider also
$$F_{K,T}(\eta_0)\stackrel{\Delta}{=} \inf \sum^{K}_{j=1}f_T(\eta_j)
\; \mbox{{\rm subject to}} \sum_{j=1}^K \eta_j \geq K \eta_0.$$
Due to relation (\ref{kappa-ext_1}), we have for
the sequence $w(r^\star_\e(b_l))$ that
 \begin{eqnarray*}
\sum_{j=1}^K \eta_j  =  \sum_{j=1}^K\E_{\theta_j}(t_{j,b_l})
  =  \frac{1}{\epsilon^2}  \sum_{j=1}^K \kappa(\theta_j,w_l)
  \geq   K (1+\delta)^2a(r_\epsilon^\star(b_l))\sim K \eta_0 \label{rela-ar-tilde-ar}.
 \end{eqnarray*}

Then, in order to obtain the same right-hand side as in (\ref{need-ext}), it is
  sufficient   to show  that for any $l$ in $\{1,\ldots,N\}$ such that $T=u_l T_d$, 
relation (\ref{pb-resoudre}) which is stated below, holds:
\begin{eqnarray}
 F_{K,T}(\eta_0) =Kf_T(\eta_0).  \label{pb-resoudre}
\end{eqnarray}
This is handled by a technical result similar to the one stated in  Lemma 7.4 and Lemma 7.5 in Ingster {\it et al. }  \cite{ITV.10}. The proof of Lemma \ref{tech-min}  is postponed to  Section \ref{subsec:appendix}.
\begin{lemma} \label{tech-min}
Set $\lambda=(T-\eta_0)f_T(\eta_0)$.
\begin{eqnarray}
\mbox{{\rm If} }  0<\eta_0 <  T -1 \mbox{ {\rm and } }  T < R < T+ ((T-\eta_0)^2 - 2\log (1+ 2 (T-\eta_0)^2))^{1/2},\label{cond-tech}
\end{eqnarray}
  then 
\begin{eqnarray}
\!\!\!\!\! \inf_{\eta \in [0,R]} (f_T(\eta) - \lambda \eta)= f_T(\eta_0) - \lambda \eta_0, \label{cond-Stat}
\end{eqnarray}
which implies that
\begin{eqnarray}
F_{K,T} (\eta_0)=  K f_T(\eta_0). \label{Y-S-Tech}
 \end{eqnarray}
\end{lemma}
As $d \rightarrow +\infty$, for any $l \in \{1,\ldots,N\}$ such that $T=u_l T_d$ with $u_l >(1 +\delta)^{2} \varphi(b) $ and  $R=  p T_d$ with $ u_l < p < u_l + \frac{u_l-(1 +\delta)^{2} \varphi(b)}{2}$, the conditions in (\ref{cond-tech}) are  then satisfied. Therefore
the application of Lemma \ref{tech-min}  yields the results since for all  ${\overline{\boldsymbol{\theta}}} \in
\Theta^{ext}_d(\tau,K^{1/2}, K^{1/2}r_\epsilon,b)$,
\begin{eqnarray*}
\E_{\overline{\boldsymbol{\theta}}} (L(u_l,b_l))
& >&  C_{u_l,b_l} K \left( \exp(-\frac{(( u_l T_d - (1+\delta)^2a(r_\epsilon^\star(b_l)))_+)^2}{2})O(1) -
\exp( -\frac{u_l T_d}{2}(1+o(1))\right), \nonumber 
\end{eqnarray*}
  which corresponds to the right-hand side of  (\ref{need-ext}).

\subsection{Lower Bound}
The prior distribution we consider is a  classical one for a functional Gaussian model. In Section 4.3 of  \cite{IS.02a} it is referred to as
the symmetric Three-point Factors. \\

\vspace{0.3cm}
\noindent {\sc Prior.}  Before defining  the prior $\Pi^{d}$ formally, we shall start with an informal discussion. \\
The prior $\Pi^{d}$ adds mass on  $(\xi_j \theta_j)_{1 \leq j \leq d}$:
the components are i.i.d. and $\xi_j$ and  $\theta_j$ are supposed to be independent.  A natural choice for $\xi_j$ is a Bernoulli with a parameter $p_d \in (0,1)$  such that  $\E (\sum_{j=1}^{d} \xi_j) \sim K$. The $\theta_j$'s are binary random variables (with probability $1/2$) such that  $\theta_j^{2}= (\theta^{\star})^{2}$ where
 the sequence $\theta^{\star}$ is a solution of the extremal problem (3.1);  this guarantees that  $\theta_j$ belongs to $\Theta(\tau, r_\epsilon)$.

Now, we define the prior distribution more precisely. Let  $\rho_d$ be any sequence of positive numbers such that
 $\rho_d \stackrel{d \rightarrow +\infty} \longrightarrow 0$ and   $d^{1-b}\rho_d^s \stackrel{d \rightarrow +\infty} \longrightarrow +\infty$, $\forall b \in (0,1)$, $\forall s >0$.
Consider  two sequences $(\xi_j)_j$ and 
$(\theta_{j,k})_{j,k}$  of independent random variables whose distributions are the following:
$$
 \left\{ \begin{array}{l}\xi_j  \sim  \mbox{ {\rm Bernoulli  $B(p_d)$ with  } }  \; p_d=d^{-b}(1+\rho_d), j \in \{1,\ldots,d\},\\
 \theta_{j,k}=\varepsilon_{j,k} \; \epsilon \;z_k,  \mbox{ {\rm with  } }   \P(\varepsilon_{j,k} =1)= \P(\varepsilon_{j,k} =-1)=\frac 1 2, \; j \in \{1,\ldots,d\}, \; k \in \BBz.
\end{array} \right. $$
The  sequence  $(z_k)_{k \in \BBz}$  is  deterministic and is defined as follows:     $(\epsilon \; z_k)_k=(\theta_k^\star)_{k \in \BBz}=\theta^\star$ where $\theta^\star$ is the sequence
that leads to
 the solution (\ref{sol_a}) of the extremal  problem  (\ref{minimisation}). In particular, this entails that
\begin{eqnarray}
\sum_{k \in \BBz} \frac{z_k^4}{2}&=& a^2(r_\epsilon) \label{a-z},  \\
(2 \pi)^{2\tau}   \sum_{k\in \BBz}| k|^{2\tau}  \;(\epsilon z_k)^2 &\leq &1, \label{sob-z}\\
\sum_{k\in \BBz}
 \; (\epsilon \;z_k)^2 &\geq &r_{\epsilon}^2 \label{sep-z}.
\end{eqnarray}
The sequences $(\xi_j)_j$ and 
$( \theta_{j,k})_{j,k}$ are also taken mutually independent.
For each $j$ in $\{1,\ldots,d\}$, we define the prior distribution $\pi_j^d$ on 
$(\xi_j,\theta_{j})$ as follows:
 \begin{eqnarray}
 \pi_j^d= ( 1-p_d ) \delta_0 + p_d \prod_{k \in \BBz} \pi_{j,k}=( 1-p_d ) \delta_0 + p_d  \pi_j, \label{level-j}
 \end{eqnarray}
 where $\pi_j=\prod_{k \in \BBz} \pi_{j,k}$, $\pi_{j,k}=   \frac 1 2  (\delta_{(-\e \;z_k)}+ \delta_{( \e \;z_k)})$ puts mass on
$ \theta_{j,k}$ 
and $\delta$ is  the Dirac mass.
Finally, we  define  the global prior $\Pi^d$  by
\begin{eqnarray*}
\Pi^d = \prod_{j=1}^d \pi^d_j \label{Prior}.
\end{eqnarray*}

\vspace{0.3cm} \noindent {\sc Minimax and Bayesian  risks.}  Denote
by $\P_{\Pi^{d}}$ the mixture of the measures $\P_{\overline{\boldsymbol{ \theta}}}$ over the
prior $\Pi^{d}$, and let $\gamma(Q)$ be the minimal
total error  probability for testing a simple null hypothesis $H_0:
\P=\P_0$ against a simple alternative $H_1: \P=Q$ regarding the measure
$\P$ of our observations $(x_{j,k})_{k\in\BBz, 1\le j\le d}$.
\begin{proposition} \label{BI}
\begin{equation}\label{Lower.1}
\gamma \ge \gamma(\P_{\Pi^d})+o(1),
\end{equation}
where $\gamma$ is the minimax total error probability
 over $\Theta_d(\tau,r_\epsilon,b)$ (see  (\ref{gamma})).
\end{proposition}

\noindent
{\it Proof of Proposition \ref{BI}.}\\
Consider
two sets $\Xi(s)$ and  $\Xi^+(s)$ defined by
$$
\Xi(s) =\{\zeta=\epsilon\,\left(\xi_1 (\varepsilon_{1,k} z_k)_k,
\ldots,  \xi_d \; (\varepsilon_{d,k} z_k)_k \right)  :\;
\sum_{j=1}^{d}\xi_j = s\},\quad 0\le s\le d,\quad
\Xi^+(s)=\bigcup_{s\le l\le d}\Xi(l).
$$
First, due to relations (\ref{sob-z}) and (\ref{sep-z}), $\Xi(K)$ is
included in $\Theta_d(\tau,r_\epsilon,b)$. This entails that
\begin{equation}\label{Lower.2}
\gamma  \geq  \gamma ( \Xi(K)).
\end{equation}
Second, let us introduce some additional priors:  for
any subset $u\subset \{1,\ldots,d\}$, define $\pi_u=\prod_{j\in
u}\pi_{j}\prod_{j\notin u}\delta_0$, where $\pi_j$ is as in (\ref{level-j}). Note that
$\pi_u$ has a support on the
collections $\zeta=\epsilon\,\left(\xi_1 (\varepsilon_{1,k} z_k)_k,
\ldots,  \xi_d \; (\varepsilon_{d,k} z_k)_k \right)$ with $\xi_j=1$
if and only if $j\in u$. For any integer $s$ such that $0\le s\le d$, let ${\cal G}_{d,s}$
be the set  of all subsets $u\subset \{1,\ldots,d\}$ of cardinality
$s$, and define $\pi_{(s)}$ as the uniform distribution on ${\cal G}_{d,s}$:
$$
\pi_{(s)}=\frac{1}{{d\choose
s}}\sum_{u\in{\cal G}_{d,s}}\pi_u.
$$
Observe that the prior $\Pi^d$ is of the form $\Pi^d=\sum_{s=0}^d
r_s\pi_{(s)}$ where $r_s=p_d^s(1-p_d)^{d-s}$.
Clearly, $\pi_{(K)}(\Xi(K))=1$, which implies
\begin{equation}\label{Lower.3}
\gamma ( \Xi(K))\ge \gamma (\P_{\pi_{(K)}}).
\end{equation}
Third, consider the conditional prior of the form $\Pi_{+}^{d}$ with respect
to $\Xi(K)^{+}$, i.e.,  $\displaystyle{\Pi_{+}^{d}(A)= \frac{\Pi^{d} (A \cap \Xi(K)^{+})}{\Pi^{d}(\Xi(K)^{+}) }}$
which is of the form
$\Pi_{+}^{d}=\displaystyle{\sum_{s=K}^d
 q_s\pi_{(s)}}$ with  $ q_s=\displaystyle{\frac{r_s}{\sum_{s=K}^{d}r_s}}$, $K \leq s \leq q$.  Let us prove that
\begin{equation}\label{Lower.4}
\gamma (\P_{\pi_{(K)}})\ge \gamma (\P_{\Pi_{+}^{d}}).
\end{equation}
Denote by ${\cal X}_K=\{(x_{j,k})_{j,k}: \;\displaystyle{\frac{d\P_{\pi_{(K)}}}{d\P_0} ((x_{j,k})_{j,k})}<1\}$ the admissible set
of the optimal test for testing $H_0: \P=\P_0$ against $H_1:
\P=\P_{\pi_{(K)}}$.
Since
$$
\g (\P_{\pi_{(K)}}) = 1 - \P_0 ({\cal X}_K) + \P_{\pi_{(K)}} ({\cal X}_K) \quad \mbox{{\rm and } }\quad  \g (\P_{\Pi_{+}^{d}}) \leq 1 - \P_0 ({\cal X}_K) + \P_{\Pi_{+}^{d}} ({\cal X}_K),
$$
proving (\ref{Lower.4}) is then reduced to checking that
\begin{equation}\label{Lower.5}
\P_{\pi_{(K)}}({\cal X}_K)\ge \P_{\Pi_{+}^{d}}({\cal X}_K).
\end{equation}
In view of Proposition 2.5 in \cite{IS.02a}, ${\cal X}_K$ is a convex
set. Also, the set ${\cal X}_K$ is sign-invariant and invariant with
respect to all permutations of the $x_{j,k}$'s; the measures
$\P_{\pi_{(s)}},\ 0\le s\le d$ have the same property of invariance with respect to  all permutations of the $x_{j,k}$'s. These
observations imply
$$
\P_{\pi_{(K)}}({\cal X}_K)=\P_{\overline{\boldsymbol{\theta}}^K}({\cal X}_K),\quad
\P_{\Pi_{+}^{d}}({\cal X}_K)=\sum_{s=K}^d q_s \P_{\overline{\boldsymbol{\theta}}^s}({\cal X}_K),
$$
where $\overline{\boldsymbol{\theta}}^s=\epsilon(\underbrace{z,\ldots,z}_{s},0,\ldots
0),\ z=(z_k)_{k\in\BBz}$. Since $\theta_{j,k}^s\ge \theta_{j,k}^K\geq 0,\
\forall \ j,k,\ s\ge K$, the application of Lemma 2.4 in \cite{IS.02a} entails that
 $\P_{\bar\theta^K}({\cal X}_K)\ge \P_{\bar\theta^s}({\cal X}_K),\
s\ge K$. This yields relation \nref{Lower.5} and hence  relation (\ref{Lower.4}).

\noindent
Finally, in view of Proposition 2.11 in \cite{IS.02a}, it remains  to check
that
\begin{eqnarray}
\g (\P_{\Pi_{+}^{d}}) =\g (\P_{\Pi^{d}}) +o(1). \label{appro_prob}
\end{eqnarray}
Similarly to the proof of Proposition 2.9. in \cite{IS.02a}, it is easily seen that  (\ref{appro_prob}) follows from the relation
\begin{eqnarray}
\Pi^d (\Xi^+(K))\stackrel{d
\rightarrow +\infty} \longrightarrow 1. \label{Masse-1}
\end{eqnarray}
Acting as in the proof of Proposition 3 in \cite{IS.02b}, we obtain
by Chebyshev's inequality,
\begin{eqnarray}
1 -\Pi^d (\Xi^+(K)) & =& \Pi^d (\sum\xi_j < d^{1-b} ) \nonumber \\
& =&  \Pi^d ( d p_d- \sum\xi_j >  d p_d-d^{1-b})   \nonumber \\
& \leq & \frac{d^{1-b}(1+\rho_d)(1-d^{-b} (1 + \rho_d))}{(d^{1-b}
\rho_d)^2}, \nonumber
\end{eqnarray}
where the ratio on the right-hand side tends to zero as
$d$  goes to infinity. Relation (\ref{Masse-1}) is then proved.

As relations \nref{Lower.2}, \nref{Lower.3}, (\ref{Lower.4}), and \nref{appro_prob}  imply \nref{Lower.1}, the proof of  Proposition \ref{BI} is completed.


\vspace{0.2cm}
\noindent
Due to Proposition \ref{BI}. the proof of the lower bound is reduced to bounding
$\gamma^\star\stackrel{\Delta}{=}\gamma ( \P_{\Pi^{d}})$ from below.

Before studying $\gamma^\star$,  we introduce some useful notation and make some helpful remarks.
Denote by  $\| \cdot \|_{TV}$  and $\| \cdot \|_2$  the distance in variation and   the $L_2$-distance between any pair of probabilities
$(P,Q)$; the latter one  is defined  by \begin{equation} \label{L2-Like}\| P-Q \|_2^2= \left\{\begin{array}{lll}  + \infty & \mbox{{\rm if }} & \mbox{{\rm $P$ does not dominate $Q$,}}\\
\E_{P} (L -1)^2&
\mbox{{\rm if }} & \mbox{{\rm $P$ dominates $Q$,}}
\end{array} \right.\end{equation}
where   $L=\displaystyle{\frac{dQ}{dP}}$ is the Radon-Nikodym derivative of $Q$ with respect to $P$.
\begin{remark} \label{Dist} Note that
\begin{itemize}
\item
 $\| \cdot \|_{TV}=\| \cdot \|_1$, where $\| \cdot \|_1$ is  the $L_1$-distance.
\item If $P$  dominates $Q$, then  $\| P-Q \|_2^2= \E_P (L^2) -1$.
\item As stated in Proposition  2.12 of \cite{IS.02a},
\begin{description}
\item If $\|P- Q \|_2=o(1)$, then $\| P - Q \|_1=o(1)$,
\item If $\|P- Q \|_2$ is bounded, then $\lim \sup \| P - Q \|_1 < 2$.
\end{description}
\end{itemize}
\end{remark}
\vspace{0.2cm}
\noindent

Using  Remark
\ref{Dist}, one has
 \begin{eqnarray}
 \mbox{ {\rm If } } \; \| \P_0 - \P_{\Pi^d}\|_{2}=o(1)  &\mbox{{\rm then }} &  \| \P_0 - \P_{\Pi^d}\|_{1} =o(1) \; \mbox{{\rm and }} \;\gamma^\star \rightarrow 1.\label{L2-LB_1}\\
\mbox{ {\rm If } } \; \| \P_0 - \P_{\Pi^d}\|_{2} = O(1) &\mbox{{\rm then }} & \lim \sup \| \P_0 - \P_{\Pi^d}\|_{1}  < 2 \; \mbox{{\rm and }} \;\lim \inf \; \gamma^\star >0. \label{L2-LB_C}
\end{eqnarray}
Therefore, if needed, the $L_2$-distance can be conveniently used
instead of the total variation distance.

Due to (\ref{L2-LB_1}) and (\ref{L2-LB_C}),
it remains to study $\| \P_0 - \P_{\Pi^d}\|_{2}$ which is  expressed
 in terms of the  Bayesian likelihood ratio $L_{\Pi^d}=\displaystyle{\frac{d\P_{\Pi^d}}{d\P_0}}$ (see relation (\ref{L2-Like})) .

\vspace{0.3cm}
\noindent
{\sc Likelihood Ratios.}
%
Here and below, when it is not absolutely necessary, we omit  the arguments of the likelihood ratios.
Then, observe  that  $L_{\Pi^d}$ is defined by:
\begin{eqnarray*}
L_{\Pi^d}  &  =&  \int \prod_{j=1}^d  (\frac{d\P_{\theta_{j}}}{d\P_0}) \; d\Pi^d \\
& =& \prod_{j=1}^d \int (\frac{d\P_{\theta_{j}}}{d\P_0}) \; d\pi_j^d  \\
& =& \prod_{j=1}^d (1 - p_d + p_d  L_j ),
 \end{eqnarray*}
where $L_j$ is the  likelihood ratio between $\P_{\pi_j}$ and $\P_0$.
Denote also by $L_{\pi_j^d}$  the likelihood ratio between $\P_{\pi_j^{d}}$ and $\P_0$ , i.e.,  $L_{\pi_j^d}=(1 - p_d + p_d   L_{j} )$.
Then   $ L_{j}$ is such that
\begin{eqnarray}
L_{j}(x_j) & =& \int \prod_{k \in \BBz}  (\frac{d\P_{\theta_{j,k}}}{d\P_0}(x_{j,k})) d\pi_j\nonumber\\
& =& \prod_{k \in \BBz} \frac 1 2 \left(  \exp(-\frac{z_k^2}{2} + z_k x_{j,k}/\e) + \exp(-\frac{z_k^2}{2} - z_k x_{j,k}/\e)\right) \nonumber\\
& =& \prod_{k \in \BBz}    \exp(-\frac{z_k^2}{2}) \cosh(z_k x_{j,k}/\e) ,\label{Lj}
\end{eqnarray}
where   $\cosh$  is  the hyperbolic cosine.
Using routine calculations, in particular, using twice the inequality $1+x \leq \exp(x), \; \forall x \in \BBr$, we
obtain
\begin{eqnarray*}
\E_0 (L_{\pi_j^d}^2) & = &   1+ p_d^2 \{\E_0 (L_j^2)-1\} \nonumber\\
& =& 1 + p_d^2 \{ \prod_{k \in \BBz} (1+2 (\sinh (\frac{z_k^2}{2}))^2)-1\} \nonumber \\
& \leq &  1 + p_d^2 \{ \exp(\sum_{k \in \BBz} 2 (\sinh (\frac{z_k^2}{2}))^2)-1\} \nonumber \\
& \leq &  \exp( p_d^2 \{ \exp(\sum_{k \in \BBz} 2 (\sinh (\frac{z_k^2}{2}))^2)-1\})
\label{RV_j},
\end{eqnarray*}
where $\sinh$ denotes the hyperbolic sine.  In view of Remark \ref{Dist}, in order to study
 $\| \P_0 - \P_{\Pi^d}\|_{2}$, it suffices to study
 $\E_0 (L_{\Pi^d}-1)^2$. The latter includes the quantity $\E_0 (L^2_{\Pi^d})$ that satisfies
\begin{eqnarray}
\E_0 (L^2_{\Pi^d}) =
\prod_{j=1}^d \; \E_0 (L^2_{\pi_j^d}) & \leq &
\exp \left(d p_d^2 \{ \exp(\sum_{k \in \BBz} 2 (\sinh (\frac{z_k^2}{2}))^2)-1\}  \right). \label{Lk}
\end{eqnarray}
\vspace{0.5cm}
As $d$ goes to infinity, the right-hand side of (\ref{Lk}) goes  to one  provided that
\begin{eqnarray}
d \; p_d^2\;\left( \exp( A) -1\right )\stackrel{d \rightarrow +\infty} \longrightarrow  0 & \mbox{{\rm with}} & A=\sum_{k \in \BBz} 2 (\sinh (\frac{z_k^2}{2}))^2. \label{RL}
\end{eqnarray}

\noindent
{\it Proof of  {\it (i)}--Theorem \ref{Moderate_case}.}\\
Recall that by assumption $b\in(0,1/2]$. We shall distinguish between two cases depending on the values of $r_\epsilon$ with respect to $r_\epsilon^\star$ defined in (\ref{sep-rates-M}). \\
{\bf Case 1:}  $r_\epsilon/r_\epsilon^\star=O(1)$. Since $dp_d^2 (a(r_{\epsilon}^\star))^{2}=O(1)$, it follows that
$dp_d^2 a^2(r_{\epsilon})=O(1)$.
Since $dp_d^2$ is bounded away from zero, $a^2(r_{\epsilon}) =O(1)$, and, due to Remark \ref{rmv0} and relations
(\ref{Lagrange}), we have $\sup_k z^2_k =o(1)$.
This entails that  $\displaystyle{\sinh^2 (\frac{z_k^2}{2}) \sim \frac{z_k^4}{4}}$, which, due to (\ref{a-z}), implies that  $A \sim \sum \frac{z_k^4}{2} \sim a^2(r_{\epsilon})$, and hence $A =O(1)$.
 It now follows  that $\exp(A) -1 \asymp A$. Finally, we  get
\begin{eqnarray}d p_d^2\;  (\exp (A) -1) \asymp d p_d^2 a^2 (r_{\epsilon})=O(1), \label{LkM} \end{eqnarray}
and the second part of  {\it (i)} in Theorem \ref{Moderate_case} is proved.

\noindent{\bf Case 2:} $r_\epsilon/r_\epsilon^\star=o(1)$.  Due to (\ref{LkM}), we have
 $\displaystyle{d p_d^2\;  (\exp (A) -1) \asymp d p_d^2 a^2 (r_{\epsilon})}$,  and since
 $\displaystyle{\frac{a^2 (r_{\epsilon})}{a^2 (r^{\star}_{\epsilon})}}=o(1)$, relation (\ref{RL}) is trivially fulfilled.

\vspace{0.5cm}

\noindent
{\it Proof of  {\it (i)}--Theorem \ref{Intermediate_case}.} \\
Now by assumption $b\in(1/2,1)$.
Due to the condition $\log (d) = o(\epsilon^{-2/(2\tau +1)})$, Remark \ref{rmv0}, and relations  (\ref{Lagrange}),
 $\sup_k  z^2_k=o(1)$. As in the moderate case, this yields
 $A \sim a^2(r_\epsilon)$, and thus
 we obtain
\begin{eqnarray}d \; p_d^2\;  (\exp (A) -1) = d p_d^2 \exp (a^2 (r_{\epsilon})(1+o(1))) \label{LkI}. \end{eqnarray}
Again, we shall consider two cases depending on the  values of $r_\epsilon$ with respect to $r_\epsilon^\star$, where $r_\epsilon^\star$
is now defined by (\ref{sep-rates-H}). \\

\noindent
{\bf Case 1:}  Suppose that $r_\epsilon/r_\epsilon^\star=o(1)$. Then   $a(r_\epsilon) =o(T_d)$. Due to equation (\ref{LkI}), this implies that
relation (\ref{RL}) is fulfilled.

\noindent
{\bf Case 2:}  Suppose that $r_\epsilon/r_\epsilon^\star=O(1)$ and
let  $c(r_\epsilon)$ be a positive constant satisfying $c^2(r_\epsilon) \log(d)=a^2(r_\epsilon)$.
Then the right-hand side of (\ref{LkI}) can be rewritten as follows:
\begin{eqnarray*}
d p_d^2 \exp (a^2 (r_{\epsilon})) & =&  d^{1-2b} (1+\rho_d)^2 \exp ( \log ( d )c^2(r_\epsilon)(1+o(1))) \\
& =& d^{1-2b +c^2(r_\epsilon)(1+o(1))}(1+\rho_d)^2. \label{LkI-Tr}
 \end{eqnarray*}
Therefore, relation (\ref{RL}) is fulfilled provided that $c(r_\epsilon)< \sqrt{2b -1}=\varphi_1(b)$, where $\varphi_1$ is defined in (\ref{phi}).
This means that a successful detection is impossible if $c(r_\epsilon) < \varphi_1(b)$, which corresponds to the intermediate sparsity case;
in fact, the  inequality $c(r_\epsilon) < \varphi_1(b)$ is valid for any $b \in (1/2, 1)$ but it
could be improved for $b \in (3/4,1)$. Indeed, for $b \in (3/4,1)$, one can show that a successful detection is impossible if
$c(r_\epsilon)$ is such that $ c(r_\epsilon) < \varphi_2(b)$, where the
function $\varphi_2$ is defined in (\ref{phi}), and for $b \in (3/4,1)$, $\varphi_1(b) < \varphi_2(b)$. That is why
the improvement is possible and  is achieved by  dealing with  a truncated version of the Bayesian  likelihood ratio $L_{\Pi^d}$.
From now, let us consider  $a(r_\epsilon) =c(r_\epsilon) \sqrt{\log d}$ with $\displaystyle{\frac{1}{\sqrt{2}}}< c(r_\epsilon) < \sqrt{2}$. The case  $c(r_\epsilon) \leq \displaystyle{\frac{1}{\sqrt{2}}}$ coincides with the intermediate sparsity case when $b\in(1/2, 3/4]$.



\vspace{0.3cm}
Thus, for some positive $v$, let us define ${\hat L}_{\Pi^d}$, the truncated likelihood ratio of $L_{\Pi^d}$:
\begin{eqnarray}
{\hat L}_{\Pi^d}= \prod_{j=1}^d {\hat L}_{\pi_j^d} = \prod_{j=1}^d  (L_{\pi_j^d}) \I_{({\tilde l_j} \leq a(r_\epsilon) \sqrt{(2+v) \log d}) }, \label{truncated-version}
\end{eqnarray}
where \begin{eqnarray}
{\tilde l_j}= \log (L_j) + \frac 1 2 \;a^2(r_\epsilon) \label{ljtilde} .\end{eqnarray}
Also put \begin{eqnarray}
l_j=\log (L_j), \label{lj}\end{eqnarray} where $L_j$ is defined by (\ref{Lj}).
Now we introduce two new probability measures  $\P_{\nu_j}$ and $\P_{\mu_j}$ expressed in terms of $\P_0$ as follows:
\begin{eqnarray}
\frac{d\P_{\nu_j}}{d\P_0}
&= &\frac{\exp(  {  l_j} )}{  \E_0(L_j)}, \label{nu}\\
\frac{d\P_{\mu_j}}{d\P_0}& = &\frac{\exp(  {2  l_j} )}{ \E_0(L_j^2)} \label{mu}.
\end{eqnarray}
In order to get a lower bound for the minimax  total error probability, it is sufficient to prove  (see the proof of Theorem 4.1 in \cite{IPT.09})  that  $\E_0(({\hat L}_{\Pi^d}-1)^2)=o(1)$,  where ${\hat L}_{\Pi^d}$ is defined in (\ref{truncated-version})  provided that
\begin{eqnarray}
\P_0 (\bigcap_{j=1}^d \; \{{\tilde l_j} \leq a(r_\epsilon) \sqrt{(2+v) \log d}\}) \rightarrow 1 . \label{proba-trunc-1}
\end{eqnarray}
In fact, it is enough to prove that
\begin{eqnarray}
 \sum_{j=1}^d \;  \P_0 ({\tilde l_j} >  a(r_\epsilon) \sqrt{(2+v) \log d}) \rightarrow 0 . \label{proba-trunc}
\end{eqnarray}
Relation (\ref{proba-trunc}), and hence relation (\ref{proba-trunc-1}),
follows from  relation (\ref{tilde_0}), which is a part of the next lemma whose  proof is postponed to
  Section \ref{subsec:appendix}. 
\begin{lemma} \label{cumulant-gj-P0} Assume that $r_\epsilon \rightarrow 0$ and 
$\log d =o(\epsilon^{-2/(2\tau +1)})$.  If $T>0$ is  such that $T=O(a^2(r_\epsilon))$,  then
\begin{eqnarray}
 \P_0 ({\tilde l_j} >  T) &\!\! \leq   & \!\! \!\! \exp \left(-\frac{T^2}{2a^2(r_\epsilon)} + o (a^2(r_\epsilon))  \right).   \label{tilde_0}\end{eqnarray}
 Moreover, if $\lim \inf (T/a^2(r_\epsilon)) >1$, then 
\begin{eqnarray}\P_{\nu_j}({\tilde l_j} > T) &\!\!  \leq & \!\! \!\! \exp \left(  -\frac{(T- a^2(r_\epsilon))^2}{2a^2(r_\epsilon)} + o (a^2(r_\epsilon))
 \right),  \label{tilde_nu} \end{eqnarray}
and if $\lim \sup (T/a^2(r_\epsilon)) <2$, then 
\begin{eqnarray}\P_{\mu_j}({\tilde l_j} \leq T)  & \!\! \leq  & \!\! \!\!\exp \left( - \frac{(T- 2a^2(r_\epsilon))^2}{2a^2(r_\epsilon)} + o (a^2(r_\epsilon))
 \right).\label{tilde_mu}
\end{eqnarray}
\end{lemma}

\noindent
Next, it remains  to  prove that $\E_0 ({\hat L}_{\Pi^d} ) \rightarrow 1$
 and  $\E_0 (({\hat L}_{\Pi^d})^2 ) \rightarrow 1$.
 This will entail the expected result that $\E_0(({\hat L}_{\Pi^d}-1)^2)=o(1)$. \\

First, consider  the term $\E_0 ( {\hat L}_{\Pi^d} )$:
\begin{eqnarray}
\E_0 ( {\hat L}_{\Pi^d} )& = & \Pi_{j=1}^d \E_0( {\hat L}_{\pi_j^d} ) \nonumber \\
& =& \Pi_{j=1}^d \E_0(1 +p_d(L_j-1) -\I_{\overline{{\cal D}_j}} \left(p_d(L_j -1)+1\right)) \nonumber \\
& =& \Pi_{j=1}^d \left(1 - p_d (\E_0 ( L_j \I_{\overline{{\cal D}_j}} )) + (-1 + p_d)\P_0 (\overline{{\cal D}_j})\right)  \nonumber \\
 & =&  \exp (\sum_{j=1}^d \log \left(1 - p_d (\E_0 ( L_j \I_{\overline{{\cal D}_j}} )) + (-1 +p_d)\P_0 (\overline{{\cal D}_j})\right), \label{hatL}
\end{eqnarray}
where  ${\cal D}_j=\{ {\tilde l_j} \leq  a(r_\epsilon) \sqrt{(2+v) \log d} \}$ and
 $\overline{{\cal D}_j} $  denotes the complement of  ${\cal D}_j$.
Relation (\ref{proba-trunc}) entails the convergence to zero of the second term in the log term of the right-hand side of (\ref{hatL}).
Therefore, in order to obtain  $\E_0 ({\hat L}_{\Pi^d} )  \rightarrow 1$, it is sufficient to prove that
\begin{eqnarray}
dp_d (\E_0 ( L_j \I_{\overline{{\cal D}_j}} )) =o(1). \label{sousP0_1}
\end{eqnarray}
Note that $\E_0 ( L_j \I_{\overline{{\cal D}_j}} )= \P_{\nu_j}(\overline{{\cal D}_j})$. Since
$\displaystyle{\frac{\sqrt{2+v}}{c(r_\epsilon)}}-1 $ is positive ($c(r_\epsilon) < \sqrt{2}$) for any positive $v$, we can applied
 relation (\ref{tilde_nu}) of Lemma \ref{cumulant-gj-P0} to get
\begin{eqnarray}
dp_d \P_{\nu_j}(\overline{{\cal D}_j})
& \leq & d p_d \exp\left(- \frac{1}{2} \log (d) ((\sqrt{2 + v} -c(r_\epsilon))^2+ o(1))\right)\nonumber \\
& =& d^{1-b}(1+\rho_d)\; d^{-\frac{1}{2}(\sqrt{2 + v} - c(r_\epsilon))^2 + o(1)},\label{fin_0}
\end{eqnarray}
where the right-hand side of (\ref{fin_0}) goes to zero as soon as $c(r_\epsilon) < \sqrt{2 + v} - \sqrt{2(1-b)}$. This yields relation (\ref{sousP0_1}).

Second, we need to study  $\E_0 ( {\hat L}^2_{\Pi^d} )$:
\begin{eqnarray}
\E_0 ( {\hat L}^2_{\Pi^d} )
& =& \prod_{j=1}^d    \E_0 ( (1 - p_d(1-L_j))^2 \I_{{\cal D}_j}) \nonumber \\
& =& \exp \left (\sum_{j=1}^d \log (1   - 2 p_d \E_0 ( (1-L_j)\I_{{\cal D}_j})+ \E_0 (p_d^2(1-L_j)^2\I_{{\cal D}_j} -\I_{\overline{{\cal D}_j}} )) \right)\nonumber .  \label{E2}
\end{eqnarray}
Since the relations $d  \P_0 (\overline{{\cal D}_j})=o(1)$  and $ d p_d \E_0 ( (1-L_j)\I_{{\cal D}_j})=o(1)$ have been  already proved,
it is sufficient  to show that
$d p_d^2 \E_0 ((1-L_j)^2\I_{{\cal D}_j})=o(1)$. 
To this end, observe that
\begin{eqnarray}
 d \E_0 (p_d^2(1-L_j)^2\I_{{\cal D}_j}) & \leq & 2  \left(d p_d^2 \P_0  ({\cal D}_j ) + d p_d^2 \E_0( L_j^2  \I_{{\cal D}_j}) \right) .\label{squared}
\end{eqnarray}
The first term on the right-hand side of (\ref{squared}) tends to zero as $d$ goes to
infinity since $d p_d^2 = d d^{-2b} (1+\rho_d)^2$ for $b\in(3/4,1)$.

To  study of the second term on the right-hand side of (\ref{squared}), we take into account the following two points: \\
(i)  since  $\sup_k  z^2_k=o(1)$, we can apply  Lemma \ref{MG_lj_0} of Section \ref{subsec:appendix} with  $h=2$, $X=x_{j,k}/\e$, and $z=z_k$, and  obtain
\begin{eqnarray}
\exp (2 \tilde{l_j}) & =& \exp (2a^{2}(r_\epsilon) ).\label{Lambda-2}
\end{eqnarray}
(ii)  since $\lim \sup (T/a^2(r_\epsilon)) < 2 $ is satisfied  as soon as $c(r_\epsilon) > \frac{\sqrt 2}{2}$
with $T=a(r_\epsilon)\sqrt{(2+v)\log d}$, we can applied relation (\ref{tilde_mu}) of Lemma \ref{cumulant-gj-P0}, which
 jointly  with  relation (\ref{Lambda-2})  leads to
\begin{eqnarray}
dp_d^2 \E_0( L_j^2  \I_{{\cal D}_j})
 & =& d d^{-2b} (1+\rho_d)^2 \P_{\mu_j}( {\tilde l}_j \leq   a (r_\epsilon) \sqrt{\log d} \;\sqrt{(2+v)}) \exp (2 \tilde{l_j} -  a^{2}(r_\epsilon) )\nonumber \\
& < & d d^{-2b} (1+\rho_d)^2  \times \nonumber \\
&& \exp \left(-\frac{a^2 (r_\epsilon)(\sqrt{2+v} \sqrt{\log d}- 2a (r_\epsilon))^2}{2 a^2 (r_\epsilon)}
+ a^2 (r_\epsilon)+ o (a^2(r_\epsilon)) \right) \nonumber \\
& =& d d^{-2b} (1+\rho_d)^2 \exp (-\frac{\log d}{2} ((\sqrt{2+v} - 2c (r_\epsilon))^2 +c^2(r_\epsilon)+o(1))) \nonumber \\
 & = &  d d^{-2b} (1+\rho_d)^2 d^{- \frac 1 2 (\sqrt{(2+v) } -2 c (r_\epsilon))^2 +c^2(r_\epsilon) +o(1)} .\label{fin_1}
\end{eqnarray}
The expression on the right-hand side of  (\ref{fin_1}) goes to zero as soon as $c(r_\epsilon) <  \sqrt{2+v } - \sqrt{2(1-b)}$.
The last inequality   is obtained by resolving the inequality
$1 - 2b - \frac 1 2 (\sqrt{2+v} -2 x)^2 +x^2  <  0$, where  $x$ is constrained to be larger than
$\frac{\sqrt 2}{2}$. This implies that a successful detection is impossible as soon as $c(r_\epsilon) < \varphi_2(b)$, where  $\varphi_2$ is defined by
(\ref{phi}).
%
%


%

 \subsection{Appendix} \label{subsec:appendix}
 \subsubsection{Proof of Lemma \ref{tech-min}.}
If there exists $\lambda$ such that (\ref{cond-Stat}) is valid,  then
  equation (\ref{Y-S-Tech}) is obtained in adapting Lemma 7.4.'s proof of  \cite{ITV.10}. Indeed, due to
 (\ref{cond-Stat}) and using the fact that $\sum_{j=1}^K \eta_j \geq K \eta_0$, we obtain
for all $\eta_j \in [0,R]$,  $ j \in \{1,\ldots,K\}$:
\begin{eqnarray}
 \sum_{j=1}^K f_T (\eta_j)  
& \geq &  \sum_{j=1}^K \inf \{  f_T (\eta_j) - \lambda \eta_j \} + \lambda K \eta_0 \nonumber \\
& \geq &  K( f_T (\eta_0) - \lambda \eta_0 )   + \lambda K \eta_0 \nonumber \\
& =& K f_T (\eta_0). \label{implique-sup}
\end{eqnarray}
On the other hand,
\begin{eqnarray}
 F_{K,T}(\eta_0)&=& \inf_{\{(\eta_1,\ldots,\eta_K): \sum \eta_j \geq K \eta_0  \}} \sum^{K}_{j=1} f_T (\eta_j) \nonumber \\
& \leq &  K f_T (\eta_0). \label{implique-inf}
\end{eqnarray}
Relations (\ref{implique-sup}) and (\ref{implique-inf}) yield relation
 (\ref{Y-S-Tech}).

\noindent
Now, let us prove that (\ref{cond-tech}) implies (\ref{cond-Stat}).
For this, set $g_T(\eta)= f_T(\eta) -\lambda \eta$ and denote by $g'_T$  and $g^{(2)}_T$ the first and second derivatives of $g_T$, respectively. Note that
   $g_T'(\eta)=(T-\eta)f_T( \eta) - \lambda $, and we  choose $\lambda=(T-\eta_0)f_T(\eta_0)$ to have
 $g'_T(\eta_0)=0$.

The study of $g_T^{(2)}$  yields that $g_T^{(2)}>0$  for  $|T-\eta| > 1$  and $g_T^{(2)} <0$  for $  |\eta-T| < 1$.
Since $0<\eta_0< T-1$, this implies that $ g'_T<0$ on $[0,\eta_0[$, $g'_T(\eta_0)=0$, $g'_T>0$ on $]\eta_0,T-1]$, $g'_T$ is decreasing on $]T-1,T+1]$,
and $g'_T$ is increasing on $]T+1,+\infty[$. Moreover, $g'_T(T-1)>0$ and  $g'_T(T)=-\lambda <0$, so that there exists $t\in ]T-1,T[$ such that $g'_T(t)=0$.
This yields  that $\eta_0$ is a local minimum of $g_T$.
In order to prove that $\eta_0$ is a global minimum of $g_T$, it is sufficient to show that $g_T(R)- g_T(\eta_0) >0$. Let us set $R=T +x$, with a positive real $x$.
We already know that
$x < T-\eta_0$ since $g_T(T + (T-\eta_0)) = f_T(\eta_0) - \lambda (T+ T-\eta_0)=f_T(\eta_0) - \lambda \eta_0 - 2 \lambda(T-\eta_0)  < g_T(\eta_0)$, where
the last inequality is valid because of the choice of $\lambda$ and $T-\eta_0$.
For  $x < (T-\eta_0)$, we  obtain
\begin{eqnarray}
g_T (R) - g_T (\eta_0) & =&
\exp(-\frac{x^2}{2}) -   (T-\eta_0)f_T(\eta_0) (T+x) -  f_T(\eta_0) +(T-\eta_0)f_T(\eta_0) \eta_0 \nonumber  \\
& > &  \exp( -\frac{x^2}{2}) - f_T(\eta_0) (2(T-\eta_0)^2+1)
 >  0\label{der},
\end{eqnarray}
where  inequality (\ref{der}) is  valid as soon as
$$\exp( -\frac{x^2}{2}) >  \exp(-\frac{(T-\eta_0)^2}{2}) (2(T-\eta_0)^2 +1) \Leftrightarrow  x < ((T-\eta_0)^2 - 2 \log (2(T-\eta_0)^2 +1) )^{1/2}.$$
Since   (\ref{cond-tech}) implies (\ref{cond-Stat}), this completes the proof of Lemma \ref{tech-min}.

 \subsubsection{Proof of Lemma \ref{cumulant-gj-P0}.}
The proof of Lemma \ref{cumulant-gj-P0} requires an additional  result stated as Lemma \ref{MG_lj_0} below.
For any $j \in \{1,\ldots,d\}$, recall that $l_j$ and ${\tilde l}_j$ are given by (\ref{lj}) and (\ref{ljtilde}), respectively.
For any $j \in \{1,\ldots,d\}$ and $k \in \BBz$, set  ${\tilde l}_{j,k} = \frac{z_k^4}{4} -\frac{z_k^2}{2} + \log (\cosh (z_k x_{j,k}/\e))$
and $l_{j,k} = -\frac{z_k^2}{2} + \log (\cosh (z_k x_{j,k}/\e))$.
Denote by $\Lambda_j $, ${\tilde \Lambda}_j$, and  ${\tilde \Lambda}_{j,k}$ the moment-generating functions of $l_j$, ${\tilde l}_j$, and ${\tilde l}_{j,k}$
under $\P_0$, respectively. 
From equations (\ref{Lj}) and (\ref{ljtilde}), it turns out that  for any $h$,
\begin{eqnarray}
{\tilde \Lambda}_j(h) &=& \prod_{k \in \BBz} {\tilde \Lambda}_{j,k}(h),  \label{Lambda-jk} \\
{\tilde \Lambda}_j(h) &=& \Lambda_j (h) \exp(h\frac{a^2(r_\epsilon)}{2}) .\label{Lambda-tilde-non}
\end{eqnarray}
Next, define the function ${\tilde g}: (z,y) \rightarrow \displaystyle{\frac{z^4}{4} - \frac{z^2}{2} +\log (\cosh(z y))}$, and
observe that the following relations hold:
\begin{eqnarray}
 {\tilde l}_{j,k} & =& {\tilde g} (z_k,x_{j,k}/\e) \label{g-tilde}, \nonumber\\
{\tilde \Lambda}_{j,k}(h) & =& \E_0 (\exp(h{\tilde g}(z_k,x_{j,k}/\e))) .\label{Generating-g-tilde}
\end{eqnarray}

\begin{lemma} \label{MG_lj_0} Let $X$ be a real standard Gaussian random variable. 
For any $z=o(1)$ and   any $h=O(1)$,
\begin{eqnarray*}
\log (\E (\exp(h{\tilde g}(z,X))))   =  h^2  \frac{z^4}{4}   + o(z^{4}).
\end{eqnarray*}
\end{lemma}

\vspace{0.3cm} \noindent
{\it Proof of Lemma \ref{MG_lj_0}.}\\
For some  $\delta>0$,
consider the event ${\cal E}=\{ |z X| < \delta\}$ and denote by $\overline{{\cal E}}$ its complement in $\BBr$.
We shall study the expectations
$G_1(h,\delta)=\E(\exp(h \log (\cosh(zX)))\I_{{\cal E}})$ and
$G_2(h,\delta)=\E(\exp(h \log (\cosh(z X)))\I_{\overline{{\cal E}}})$ separately.
At this point, we choose $\delta$ small enough ($\delta =o(1)$) to satisfy $z \delta^{-1}=o(1)$.

First, let us  study the term $G_2(h,\delta)$. With the use of the  inequality $\cosh(x) \leq \exp(|x|), \; \forall x \in \BBr, $ and the fact that $h=O(1)$, 
the routine calculations of exponential moments of a real Gaussian random variable lead to
\begin{eqnarray}
G_2 (h,\delta)& \leq & \E(\exp (h|zX|)\I_{(|zX| \geq \delta)})  \nonumber \\
& =& 2 \E (  \exp (h zX )\I_{(X \geq \delta/z)}) \nonumber \\
& =& \frac{2}{\sqrt{2 \pi}} \int_{\BBr^{+}} \exp (-\frac 1 2 (x-hz)^2) \; \I_{(x \geq \frac{\delta}{z})} dx \exp (\frac 1 2 h^2 z^2)  \nonumber \\
& \leq & 2 \exp(h^2 \frac{z^2}{2}) \exp( -\frac 1 2 (\frac{\delta}{z} -hz)^2 )  \nonumber \\
& \leq  & 2 \exp\left( -\frac 1 2 \frac{\delta^2}{z^2} + o(1)\right),\label{G2}
\end{eqnarray}
where, with our choice of $\delta$,  the right-hand side of (\ref{G2}) is  small. 

Now, we move on to the term $G_1(h,\delta)$. If $\delta$  is small enough, then      $|zX|$ is also small  and the following relation holds:
\begin{eqnarray}
 \log (\cosh (zX))& =& 
\frac{z^2}{2} X^2   -\frac{z^4 }{12} X^4 + o(z^4 X^4). \label{Dv_lj}
\end{eqnarray}
Then the routine calculations of exponential moments as above lead to the following:
\begin{eqnarray}
 G_1(h,\delta)&=&  \E\left(\exp\left(h (\frac{z^2}{2} X^2   -\frac{z^4 }{12} X^4(1+o(1)) )\right) \I_{{\cal E}}\right)\nonumber \\
 & =&\E\left(\exp(h (\frac{z^2}{2} X^2 )) (1  -h\frac{z^4 }{12} X^4  (1+o(1)))  \I_{{\cal E}} \right)\nonumber \\
 &=&  \exp(-\frac 1 2 \log (1 -h z^2) )  \exp( - \frac{h }{4}z^4 (1+ o (1)))\nonumber \\
 & =& \exp (  \frac h 2 z^2  + \frac{ h^2}{ 4}  z^4 (1+o(1))) \exp( - \frac{h }{4}z^4 (1+ o (1))) \nonumber \\
 & = &   \exp (  \frac h 2 z^2  + \frac{h^2}{ 4}  z^4 - \frac{h }{4}z^4  + o( z^4  )).  \label{G1}
\end{eqnarray}

Taking
 $h=O(1)$,    $z=o(1)$, $\delta=o(1)$ and $z \delta^{-1}=o(1)$ in  relations (\ref{G2}) and (\ref{G1}) entails that
$G_1(h,\delta)=O(1)$,  $G_2(h,\delta)=O(\exp(-\delta^2/(2z^2))=o(1)$, and therefore $G_2(h,\delta)(G_1(h,\delta))^{-1}=o(1)$.

Next, due to (\ref{G2}), (\ref{G1}) and using the fact that $h=O(1)$, $z=o(1)$,  for small $\delta$  such that $z_0 \delta^{-1}=o(1)$ and $\delta =o(1)$, we obtain
 \begin{eqnarray}
\log (\E (\exp(h{\tilde g}(z,X))))  & = &   \log (   G_1(h,\delta) + G_2(h,\delta))  -\frac h 2 (z^2 -\frac{z^{4}}{2})  \nonumber \\
& =&     (\log G_1(h,\delta) -\frac h 2 (z^2 -\frac{z^{4}}{2}) ) + \log  ( 1 + \frac{G_2(h,\delta)}{G_1(h,\delta) }) \nonumber \\
& =&    h^2 \frac{z^4}{4}   + o(z^{4}) +  \frac{G_2(h,\delta)}{G_1(h,\delta) }(1 +o(1)) \nonumber \\
& =&     (h^2 \frac{z^4}{4}   +  o(z^{4}) )(1 +  \frac{G_2(h,\delta)(1 +o(1))}{G_1(h,\delta)(h^2 \frac{z^4}{4}   + o(z^{4})  )}  )\nonumber \\
  & =&  h^2 \frac{z^4}{4} +  o(z^{4}), \label{mgf}
 \end{eqnarray}
where  relation (\ref{mgf}) holds  provided that  \begin{eqnarray}
 \frac{G_2(h,\delta)}{G_1(h,\delta) (h^2 \frac{z^4}{4} + o(z^{4}))  }=o(1) \label{tech2}. \end{eqnarray}
It is then sufficient to  prove (\ref{tech2}) since (\ref{mgf}) is the expected result of Lemma \ref{MG_lj_0}.
Recall that  $h=O(1)$ and $z=o(1)$ entail that    $G_1(h,\delta)=O(1)$ and $G_2(h,\delta)=O(\exp(-\delta^2/(2 z^2)))$. Then, it is sufficient
to establish that
$\displaystyle{\exp( -\frac 1 2 \frac{\delta^2}{z^2} ) z^{-4}}=o(1)$. The latter holds if we choose $\delta$   such that
$ \delta^{-1}=o((z \sqrt{\log(z^{-1})})^{-1})$.

\vspace{0.5cm}
\noindent
{\it Proof of Lemma \ref{cumulant-gj-P0}.} \\
Remark \ref{rmv0} and relations (\ref{Lagrange}) imply that  $\displaystyle{\sup_k} z_k^2 \leq z_0^2= o(1)$ as soon as
$\log (d) = o (\epsilon^{-2/(2\tau +1)})$. Due to (\ref{Generating-g-tilde}), for any  $h$ such that $h=O(1)$,
 Lemma \ref{MG_lj_0} can be applied to the moment-generating function
$\Lambda_{j,k}(h)$.

Here and later, we consider any $j \in \{1,\ldots,d\}$ and any $k \in \BBz$.
Due to relations (\ref{Lambda-jk}),  (\ref{Generating-g-tilde}),  (\ref{a-z}), (\ref{ljtilde}), by applying Lemma \ref{MG_lj_0}  and using the exponential Chebyshev's inequality,
we obtain for any positive  $h$ such that $h =O(1)$,
\begin{eqnarray}
\P_0 ({\tilde l_j} >  T) & \leq & {\tilde \Lambda}_j(h) \exp (-hT)  \nonumber \\
& \leq & \exp( \frac{h^2}{2} a^2(r_\epsilon) - hT  +o (a^2(r_\epsilon)) ). \label{Min_0}
\end{eqnarray}
The minimum on the right-hand side of (\ref{Min_0}) is  attained for
 $h=\displaystyle{\frac{T}{a^2(r_\epsilon)}}$  which is positive and of order $1$; this allows us to prove relation (\ref{tilde_0}).

Due to  relations  (\ref{Generating-g-tilde}), (\ref{Lambda-jk}), (\ref{a-z}), (\ref{ljtilde}), (\ref{nu}), by applying again Lemma \ref{MG_lj_0} and
using the exponential Chebyshev's inequality,  we obtain  for any positive  $h$ such that $h =O(1)$,
\begin{eqnarray}
\P_{\nu_j} ({\tilde l_j} >  T) & \leq & \E_{\nu_j} (\exp ({\tilde l}_j h)) \exp (-hT) \nonumber \\
& =& {\tilde \Lambda}_j(h+1) \exp ( -\frac{a^2(r_\epsilon)}{2}-hT) \nonumber \\
& =& \exp\left(  \frac{(h+1)^2}{2}a^2(r_\epsilon)-\frac{a^2(r_\epsilon)}{2} -hT  + o (a^2(r_\epsilon)) \right), \label{Min_1}
\end{eqnarray}
where the minimum in the right-hand side of (\ref{Min_1}) is attained for  $h=\displaystyle{\frac{T}{a^2(r_\epsilon)}}-1$ which is  positive and of order $1$; this yields relation (\ref{tilde_nu}).

Recall that under the assumption of Lemma \ref{cumulant-gj-P0}, the quantity $2a^2(r_\epsilon)-T$ is positive.
Therefore, from (\ref{Generating-g-tilde}),    (\ref{Lambda-jk}), (\ref{a-z}), (\ref{mu}), (\ref{ljtilde}), and (\ref{Lambda-tilde-non}), 
applying Lemma \ref{MG_lj_0}  and using
the exponential Chebyshev's inequality, we get for any positive  $h$ such that $h =O(1)$,
\begin{eqnarray}
\P_{\mu_j} ({\tilde l_j} \leq  T) &=& \P_{\mu_j} (-{\tilde l_j} \geq -T) \nonumber \\
&=& \E_{\mu_j} (\exp ( -{\tilde l}_j h)) \exp (hT)) \nonumber \\
&=& \E_{0} (\exp ( -{\tilde l}_j h) \exp(2 {\tilde l}_j)) \exp(-a^2(r_\epsilon)) (\Lambda_j(2))^{-1}\exp (hT) \nonumber \\
& =& {\tilde \Lambda}_j(2-h)  (\Lambda_j(2))^{-1} \exp ( - a^2 (r_\epsilon) + T h) \nonumber \\
& =& {\tilde \Lambda}_j(2-h)  ({\tilde \Lambda}_j(2))^{-1} \exp(a^2 (r_\epsilon)) \exp ( - a^2 (r_\epsilon) + T h) \nonumber \\
& =& \exp\left( \frac 1 2  (2 - h)^2 a^2(r_\epsilon) -2a^2(r_\epsilon)  + T h + o (a^2(r_\epsilon))  \right),  \label{Min_2}
\end{eqnarray}
where the minimum in the right-hand side of (\ref{Min_2}) is  achieved for  $h=-\displaystyle{\frac{T}{a^2(r_\epsilon)}}+2$ which is  positive and  of order $O(1)$;
this yields relation (\ref{tilde_mu}). The proof of Lemma \ref{cumulant-gj-P0} is completed.

\vspace{1cm}
\noindent
\textbf{Acknowledgements:}  We thank one of our colleagues for his contribution to improve the English language of the paper.
%
%
%
%
%

\end{document}